\documentclass[11pt,reqno,a4paper]{article}
\usepackage{geometry}
\usepackage{amsmath,amsthm,amstext,amscd,amssymb,euscript,url}
\usepackage{mathrsfs}
\usepackage{calrsfs}
\usepackage{epsfig}
\usepackage[inline]{enumitem}
\usepackage{todonotes}
\usepackage{mathtools}
\usepackage[absolute,overlay]{textpos}
\usepackage{graphicx}
\usepackage{subfig}
\usepackage{tikz}
\usetikzlibrary{matrix}
\usetikzlibrary{arrows,positioning}

\usepackage{color}

\newcommand{\gi}{g}

\newcommand{\pii}{\Pi}

\renewcommand{\P}{\mathbb{P}}

\newcommand{\Z}{\mathbb Z}
\newcommand{\R}{\mathbb R}
\newcommand{\N}{\mathbb N}

\newcommand{\E}{\mathbb E}
\newcommand{\Zd}{\mathbb Z^d}
\newcommand{\Rd}{\mathbb R^d}

\renewcommand{\Pr}{\mathbb P}

\renewcommand{\phi}{\varphi}

\newcommand{\La}{\ensuremath{\Lambda}}
\newcommand{\la}{\ensuremath{\lambda}}
\newcommand{\si}{\ensuremath{\sigma}}

\newcommand{\loc}{\mathcal{L}}

\newcommand{\Di}{{\mathfrak D}}
\newcommand{\di}{{\mathfrak d}}
\newcommand{\biy}{{\mathbf y}}

\newcommand{\cs}{c_\sigma}

\def\1{{\mathchoice {\rm 1\mskip-4mu l} {\rm 1\mskip-4mu l}
{\rm 1\mskip-4.5mu l} {\rm 1\mskip-5mu l}}}

\newtheorem{theorem}{{\small T}{\scriptsize HEOREM}}[section]
\newtheorem{corollary}{{\bf{\small C}{\scriptsize OROLLARY}}}[section]
\newtheorem{proposition}{{\bf{\small P}{\scriptsize ROPOSITION}}}[section]
\newtheorem{lemma}{{\bf{\small L}{\scriptsize EMMA}}}[section]
\newtheorem{remark}{{\bf{\small R}{\scriptsize EMARK}}}[section]
\newtheorem{definition}{{\bf{\small D}{\scriptsize EFINITION}}}[section]
\newtheorem{induction}{{\bf{\small I}{\scriptsize NDUCTIVE HYPOTHESIS}}}[section]

\renewenvironment{proof}[1]
{\noindent{{\bf{\small{ P}{\scriptsize ROOF}}}.}\hspace{0.1cm} #1} {$\;\qed$\newline}

\newcommand{\beq}{\begin{eqnarray}}
\newcommand{\eeq}{\end{eqnarray}}

\newcommand{\ba}{\begin{align*}}
\newcommand{\ea}{\end{align*}}

\newcommand{\be}{\begin{equation}}
\newcommand{\ee}{\end{equation}}

\newcommand{\bl}{\begin{lemma}}
\newcommand{\el}{\end{lemma}}

\newcommand{\br}{\begin{remark}}
\newcommand{\er}{\end{remark}}

\newcommand{\bt}{\begin{theorem}}
\newcommand{\et}{\end{theorem}}

\newcommand{\bd}{\begin{definition}}
\newcommand{\ed}{\end{definition}}

\newcommand{\bind}{\begin{induction}}
\newcommand{\eind}{\end{induction}}

\newcommand{\bp}{\begin{proposition}}
\newcommand{\ep}{\end{proposition}}

\newcommand{\bc}{\begin{corollary}}
\newcommand{\ec}{\end{corollary}}

\newcommand{\bpr}{\begin{proof}}
\newcommand{\epr}{\end{proof}}

\newcommand{\bi}{\begin{itemize}}
\newcommand{\ei}{\end{itemize}}

\newcommand{\ben}{\begin{enumerate}}
\newcommand{\een}{\end{enumerate}}

\newcommand{\caE}{{\mathrsfs E}}

\newcommand{\caF}{{\mathcal F}}
\newcommand{\caG}{{\mathcal G}}

\newcommand{\caP}{{\mathcal P}}

\newcommand{\caR}{{\mathcal R}}

\newcommand{\caW}{{\mathcal W}}
\newcommand{\caY}{{\mathcal X}}
\newcommand{\caX}{{\mathcal Y}}

\newcommand{\IRW}{\text{\normalfont IRW}}
\newcommand{\SIP}{\text{\normalfont SIP}}

\newcommand{\SEP}{\text{\normalfont SEP}}

\newcommand{\bix}{{\mathbf x}}

\renewcommand{\(}{\left(}
\renewcommand{\)}{\right)}
\newcommand{\nn}{\nonumber}

\newmuskip\pFqmuskip

\newcommand\pFq[6][8]{%
	\begingroup % only local assignments
	\pFqmuskip=#1mu\relax
	\mathcode`\,=\string"8000
	\begingroup\lccode`\~=`\,
	\lowercase{\endgroup\let~}\pFqcomma
	{}_{#2}F_{#3}{\left[\genfrac..{0pt}{}{#4}{#5};#6\right]}%
	\endgroup
}
\newcommand{\pFqcomma}{\mskip\pFqmuskip}

\newcommand*{\myprime}{^{\prime}\mkern-1.2mu}

\newcommand{\norm}[1]{\left\lVert#1\right\rVert}

\newtheorem{asu}{Assumption}
\newcommand{\basu}{\begin{asu}}
\newcommand{\easu}{\end{asu}}

\begin{document}
\title{{\bf Higher order fluctuation fields and orthogonal duality polynomials}}
\author{Mario Ayala, Gioia Carinci and Frank Redig
% Wioletta Ruszel$^{\textup{{\tiny(a)}}}$,
\\
%\small $^{\textup{(a)}}$
\small{Delft Institute of Applied Mathematics}\\
\small{Delft University of Technology}\\
{\small Van Mourik Broekmanweg 6, 2628 XE Delft}
\\
\small{The Netherlands}
%\\
}
\maketitle

\begin{abstract}
Inspired by the works in \cite{assing2007limit} and \cite{gonccalves2019quadratic} we introduce what we call $k$-th-order fluctuation fields and study their scaling limits. This construction is done in the context of particle systems with the property of orthogonal self-duality. This type of duality provides us with a setting in which we are able to interpret these fields as some type of discrete analogue of powers of the well-known density fluctuation field. We show that the weak limit of the $k$-th order field satisfies a recursive martingale problem that formally corresponds to the SPDE associated with the $k$th-power of a generalized Ornstein-Uhlenbeck process.
\end{abstract}

%\newpage
\tableofcontents

%\newpage
\section{Introduction}
In the context of interacting particle systems with a conserved quantity (such as the number of particles) in \cite{de1984survey, kipnis2013scaling} one studies the time-dependent density fluctuation field
\[
\caY^{(n)}(\phi, \eta(n^2 t)) =  \frac{1}{n^{d/2}} \sum_{x \in \Zd} \phi(x/n) (\eta_x (n^2 t)- \rho).
\]
Here $\phi$ denotes a test-function, and $\eta_x$ the number of particles at site $x\in\Zd$. 
The quantity $\caY^{(n)}(\phi, \eta(n^2 t))$ is then considered as a random time-dependent (Schwartz) distribution. In
 a variety of models with particle number conservation (such as zero  range processes, simple exclusion process, etc.),  this time-dependent field is proved to converge, at equilibrium, to a stationary infinite dimensional Ornstein-Uhlenbeck process. 
 This scaling limit behavior of the density fluctuation field can be thought of as a generalized central limit theorem, and, as such, as a correction or refinement of  the hydrodynamic limit, which can be thought of as a law of large numbers.\\

The usual strategy of proof (see e.g. Chapter 11 of \cite{kipnis2013scaling}) is to start from the Dynkin martingale associated to the density field and prove convergence of the drift term via the Boltzmann-Gibbs principle (the drift term becomes in the scaling limit a function of the density field), and convergence of the noise term via characterization of its quadratic variation (which becomes deterministic in the scaling limit). This then eventually leads to the informally written SPDE
\[
d \caY_t= D  \Delta \caY_t +   \si(\rho) \nabla d \caW_t
\]
where $\rho$ is the parameter of the invariant measure associated to the density, $\Delta$ denotes the Laplacian, and where $\si(\rho) \nabla d \caW_t$ is an informal notation for Gaussian white noise with variance $\si^2(\rho) \int (\nabla \phi)^2 dx$.\\

In interacting particle systems with (self-)duality, the drift term in the equation for the density field is already microscopically (i.e., without rescaling) a (linear) function of the density field.  As a consequence, closing the equation and proving  convergence to the limiting Ornstein-Uhlenbeck process, is, for self-dual systems, particularly simple  and don't require the use of a Boltzmann-Gibbs principle. This simplification suggests that, in that context, we can obtain more detailed results about fluctuation fields of more general observables.
Orthogonal polynomial duality is a useful tool in the study of fluctuation fields, and associated Boltzmann-Gibbs principles, as we have seen in  \cite{ayala2018quantitative}.

The density fluctuation field can be viewed as the lowest (i.e., first) order of a sequence of fields associated to orthogonal polynomials. Indeed, in all the models with orthogonal polynomial self-duality, the function $\( \eta_x-\rho \)$ is the first order orthogonal polynomial up to a multiplicative constant. Orthogonal polynomials are indexed by finite particle configurations, i.e.,  the dual configurations. If we denote by $D(x_1,  \ldots x_k; \eta)$ the orthogonal polynomial associated to the dual configuration $ \sum_{i=1}^n  \delta_{x_i}$, then a natural field generalizing the density fluctuation field is
\[
\caY^{(n,k)}(\Phi,\eta) =n^{-kd/2} \sum_{x_i \in \Zd} D(x_1,\ldots, x_k; \eta) \cdot \Phi\left(\tfrac{x_1}{n},\ldots, \tfrac{x_k}{n}\right).
\]
In the context of  exclusion process the case $k=2$ (orthogonal polynomial of order 2) has been studied in  \cite{gonccalves2019quadratic}, where this field,  called the quadratic fluctuation field,   is shown to converge, in the limit $n \to \infty$, to the solution of a martingale problem.  The quadratic variation  of this $2$nd-order field is proven to be a function of the $1$st order field (the  density field). 
 From the result on the quadratic (k=2) field one can conjecture the existence of a more general structure where the $k$th-order orthogonal polynomials field  satisfies, in the scaling limit, a martingale problem with quadratic variation depending on the $k-1$-order  field.

In this paper we show exactly the emergence of a scenario of this type: within a general class of models with orthogonal polynomial self-duality we consider the fluctuation fields associated to orthogonal polynomials and prove that they converge, in the scaling limit, to the solution of a recursive system of martingale problems. We believe that this can also be a first step in the direction of defining non-linear fields, such as the square of the density field, via approximation of the identity, i.e., via a singular linear observable (cf. \cite{gonccalves2019quadratic}) of the field constructed in our paper.\\

The rest of our paper is organized as follows. In Section 2 we define the basic models, and introduce orthogonal polynomial duality. In Section 3 we define the fluctuation fields, in Section 4 we introduce a coordinate version of the dual process, a  technical tool that will prove to be useful  later on. In Section 5 we state the main result, Theorem \ref{mainkthorder} below,  and outline a strategy of its inductive proof. Finally, the rest of the sections are devoted to the proof of Theorem \ref{mainkthorder}.

\section{The models}

\subsection{The infinite configuration process} \label{inficonf}
We consider an interacting particle system where an infinite number of particles randomly hop on the lattice $\Zd$.
Configurations are denoted by $\eta, \xi, \zeta$ and are elements of $\Omega\subseteq \N^{\Zd}$ (where
$\N$ denotes the natural numbers including zero).
We denote by $\eta_x$ the number of particles at $x$ in the configuration $\eta\in\Omega$.
We have in mind  symmetric processes of the type independent random walkers,  inclusion or exclusion.
We fix two parameters $(\sigma,\alpha)\in \{0,1\}\times [0,\infty)\cup \{-1\}\times \N$ and we define the  generator working on local functions $f:\Omega\to\R$ as
\be\label{gen}
\loc f(\eta) = \sum_{i \in \Zd} \sum_{r \in \Zd} p(r) \eta_i ( \alpha + \sigma \eta_{i+r})(f(\eta^{i,i+r})-f(\eta))
\ee
where $\eta^{i,i+r}$ denotes the configuration obtained from $\eta$ by removing a particle from $i$ and putting it at $i+r$.
The state space $\Omega$ has to be defined and its form depends on the choice of the parameters $\alpha$ and $\sigma$.
\noindent
We assume that $p(r)$ is a symmetric, finite range, irreducible Markov transition function on $\Zd$:
\ben
\item \underline {Symmetry}. The function $p:\R^d\to [0,\infty)$ is  of the form:
\beq\label{ps}
p(r_1,\ldots, r_d)= p(|r_1|, \ldots, |r_d|)
\eeq
and such that $p(r_{\sigma(1)}, \ldots, r_{\sigma(d)})=p(r_1, \ldots, r_d)$ for all $\sigma\in {\cal P}(d)$, the set of permutations of $\{1, \ldots, d\}$.

\item \underline {Finite-range}. There exists a finite subset $\caR \subset \Zd$  of the form $\caR=[-R,R]^d\cap \Zd$, for some $R\in \N$, $ R>1$, such that $p(r)=0$ for all $r \notin \caR$.

\item \underline {Irreducibility}. For all $x,y\in\Zd$ there exists a sequence $i_1=x,\ldots,i_n=y$ such
that 
\[
\prod_{k=1}^{n-1} p(i_k-i_{k+1})>0.
\]
\een
We will also assume, without loss of generality, that $p(0)=0$, and denote by $\chi$ the second moment:
\be\label{chidef}
\chi := \sum_{r \in \caR} r_\ell^2 \cdot p(r), \qquad \text{for all} \quad \ell\in \{1,\ldots,d\}.
\ee
\vskip.2cm
\noindent
For the associated Markov processes on $\Omega$, we use the notation $\{\eta(t):t\geq 0\}$, 
$\eta_x(t)$ denoting the number of particles at time $t$ at location $x\in \Zd$.
These particle systems have a one-parameter family of homogeneous (w.r.t. translations) reversible and ergodic product measures $\nu_\rho, \rho>0$, indexed by the density of particles, i.e.,
\be
\int\eta_0 d\nu_\rho= \rho.
\ee
The nature of the underlying dynamics and the type of reversible measure we obtain is regulated by the parameter $\sigma \in \Z$ as follows.
\begin{description}
\item[Independent random walkers ($\IRW$):] This particle system corresponds to the choice $\sigma=0$ and the intensity parameter $\alpha \in \R $ regulates the rate at which the  walkers move. The reversible measures $\nu_{\rho}$, $\rho>0$ are products of Poisson distributions with parameter $\rho$, $\nu_\rho =\otimes_{i\in \Zd} \text{Pois}(\rho)$, i.e. the marginals are given by
\[
\mathbb P_{\nu_\rho}(\eta_i=n)= \frac{1}{Z_{\rho}} \cdot\frac{\rho^n}{n!}, \qquad \qquad  Z_{\rho}= e^{-\rho},
\qquad \forall \: i\in \Zd.
\]

\item[Symmetric exclusion process ($\SEP(\alpha)$):] The choice $\sigma=-1$ results in exclusion interaction. For this process the  parameter $\alpha$ takes values  in the set of natural number, $\alpha \in \N$, as it determines the maximum number of particles allowed  per site. This system is well known to have reversible measures $\nu_\rho$, $\rho\in(0,\alpha)$, that are products of  Binomial distributions: $\nu_\rho=\otimes_{i\in \Zd} \text{Binom}\(\alpha,\tfrac \rho \alpha\)$ whose  marginals are given by 
\[
\mathbb P_{\nu_\rho}
(\eta_i=n)= \frac{1}{Z_{\alpha,\rho}} \cdot \binom{\alpha}{n} \cdot \left(\frac{\rho}{\alpha- \rho} \right)^{n}, \qquad  \qquad Z_{\alpha,\rho}= \(\frac \alpha{ \alpha-\rho}\)^{\alpha},
 \qquad \forall \: i\in \Zd.
\]

\item[Symmetric inclusion process ($\SIP(\alpha)$):]   The choice $\sigma=1$ gives rise to an interaction of inclusion-type consisting of  particles  attracting each other. The SIP is known to have  products of  Negative-Binomial distributions as reversible measures, i.e. $\nu_\rho$, $\rho>0$ with $\nu_\rho=\otimes_{i\in \Zd} \text{Neg-Binom}\(\alpha,\tfrac{\rho}{\rho+\alpha}\)$ with marginals
\[
\mathbb P_{\nu_\rho}(\eta_i=n)= \frac{1}{Z_{\alpha,\rho}} \cdot  \frac{\Gamma(\alpha +  n)}{\Gamma(\alpha) \cdot n!}  \(\frac{\rho}{\alpha +\rho}\)^{n}, \qquad \qquad Z_{\alpha,\rho}= \(\frac{\alpha+\rho}{\alpha}\)^{\alpha},
 \qquad \forall \: i\in \Zd.
\]

\end{description}

\vskip.2cm
\noindent
The definition of the state space $\Omega$ is different on each case, depending on whether there are restrictions or not on the total number of particles allowed per site.
This is finite  for the exclusion process, thus, for SEP$(\alpha)$, we have $\Omega=\{0, 1, \ldots, \alpha\}^{\Zd}$. 
The situation is different in the cases of IRW and  SIP,  for which, in principle, there are no restrictions. Nevertheless, one has to avoid explosions of the number of particles  in a given site.  For this reason the characterization of $\Omega$ in these cases (i.e. for $\sigma\ge 0$) is a more subtle problem whose treatment is beyond the scope of this paper. Here we will restrict to   implicitly define $\Omega$ as the set  of   configurations in $\N^{\Zd}$ whose evolution  $\eta(t)$ is well-defined and belonging to $\Omega$ for all subsequent  times $t\ge 0$.
For what we are concerned, a precise characterization of $\Omega$ is not necessary as
we will  mostly work in $L^2(\nu_\rho)$ spaces,
and the set $\Omega$ (implicitly defined above) has $\nu_\rho$ measure $1$ for all $\rho$. 
%\col{Conditions on $p(\cdot)$ in order to avoid escape to $\infty$?}

\subsection{The finite configuration processes}
The process introduced in Section \ref{inficonf} can also be realized with a fixed finite number of particles. For a process with $k\in \N$ particles we denote by $\Omega_k$ its state space, more precisely:
\be
\Omega_k = \Big\{ \xi \in \Omega : \|\xi\|:= \sum_{x \in \Zd} \xi_x = k \Big\}.
\ee
We will then denote by $\{ \xi(t) : t \geq 0 \}$ the $\Omega_k$-valued Markov process, with infinitesimal generator given by
\be\label{genk}
\loc^{(k)} f(\xi) = \sum_{i \in \Zd} \sum_{r \in \caR} p(r) \xi_i ( \alpha + \sigma \xi_{i+r})(f(\xi^{i,i+r})-f(\xi))
\ee
working on functions $f:\Omega_k\to \R$.\\
\noindent
We now define the measure $\La(\cdot)$, which is a product measure not depending on $k$:
\be\label{muk}
\La(\xi)=\prod_{i\in \Zd} \la(\xi_i) \qquad \text{with} \quad \la(m)=
\left\{
\begin{array}{llll}
 \frac 1{m!}, &m\in \N& \text{for }\sigma=0& \text{IRW}\\
 \\
 \frac{\alpha!}{m!(\alpha-m)!}, & m \in \{0,\ldots, \alpha\} & \text{for }\sigma=-1& \text{SEP}(\alpha)\\
 &\\
\frac{m!\Gamma(\alpha + m)}{\Gamma(\alpha)},& m\in \N &\text{for }\sigma=1 & \text{SIP}(\alpha)
\end{array}
\right.\ee
Notice that by detailed balance we can verify that the measure $\La(\cdot)$ is reversible. As a consequence of this reversibility  we can infer that the $k$-particles generator $\loc^{(k)}$ is self-adjoint with respect to the inner product $\langle\cdot,\cdot\rangle_{\La}$, i.e. for all $f, g \in {L}^2(\Omega_k, \La)$ we have
\be
\langle f, \loc^{(k)} g \rangle_{\La} = \langle \loc^{(k)}f, g \rangle_{\La}
\ee
where the inner product is given by
\be
\langle f,g \rangle_{\La} = \sum_{\xi \in \Omega_k} f(\xi) g(\xi) \La(\xi).
\ee

\subsection{Orthogonal polynomial self-duality}
The processes defined in Section \ref{inficonf} share a self-duality property that will be crucial in our analysis.  Define the set
\be
\Omega_f = \mathop{\bigcup}_{ k \in \N} \Omega_k
\ee
 of configurations with a finite number of particles,  the self-duality functions that we consider in this paper   are functions $D_\rho:\Omega_f\times \Omega \to\R$ parametrized by the density $\rho>0$  satisfying the  following properties.
\ben
\item \underline {Self-duality}:
\be\label{expduality}
\E_\eta \left[D_\rho(\xi,\eta(t))\right]=\E_\xi \left[D_\rho(\xi(t), \cdot)\right] \qquad \text{for all} \quad  \xi\in \Omega_f, \eta\in \Omega
\ee
or, equivalently,  
\be\label{Genduality}
[\loc D_\rho(\xi,\cdot)](\eta) =[\loc^{(k)} D_\rho(\cdot, \eta)](\xi) \qquad \text{for all} \quad  \xi\in \Omega_f, \eta\in \Omega.
\ee
\item \underline{Factorized polynomials}:
\[
D_\rho(\xi,\eta)=\prod_{i\in\Zd} d_\rho(\xi_i,\eta_i)
\]
where $d_\rho(0,n)=1$, and $d_\rho(k,\cdot)$ is a polynomial of degree $k$.
\item \underline{Orthogonality}:
\be\label{Ort}
\int D_\rho(\xi, \eta) D_\rho(\xi', \eta) \, d\nu_\rho(\eta) =  \delta_{\xi,\xi'} \cdot \frac{1}{\mu_\rho(\xi)}
\ee
where
\be\label{defmurho}
\mu_{\rho}(\xi) := \( \int D_\rho(\xi, \eta)^2 \, d\nu_\rho(\eta) \)^{-1}.
\ee
\een
For proofs of self-duality with these functions we refer to \cite{franceschini2017stochastic} and \cite{redig2018factorized}.
\br
Notice that, as a consequence of the orthogonality property \eqref{Ort}, we have that
\beq\label{ort}
\int \E_{\eta}\left[D_\rho(\xi,\eta(t))\right]\cdot D_\rho(\xi',\eta) \; d\nu_\rho(\eta)=p_t(\xi,\xi') \cdot \frac{1}{\mu_\rho(\xi)}
\eeq
where $p_t(\cdot,\cdot)$ is the transition probability function of the dual process $\{\xi(t): t\ge 0\}$. Moreover, if we use the reversibility of the measure $\nu_\rho$ on the LHS of \eqref{ort} we obtain
\beq\label{ort2}
p_t(\xi,\xi') \cdot \frac{1}{\mu_\rho(\xi)} &=& \int \E_{\eta}\left[D_\rho(\xi,\eta(t))\right]\cdot D_\rho(\xi',\eta) \; d\nu_\rho(\eta) \nn \\
&=& \int D_\rho(\xi,\eta) \cdot \E_{\eta}\left[ D_\rho(\xi',\eta(t)) \right] \; d\nu_\rho(\eta) \nn \\
&=& p_t(\xi',\xi) \cdot \frac{1}{\mu_\rho(\xi')}
\eeq
which, by detailed balance, implies the reversibility of the measure $\mu_\rho(\xi)$. This in turn implies that there exists a constant $c(k,\rho)$ such that
\beq\label{La}
\La(\xi)=c(k,\rho)\cdot \mu_\rho(\xi) \qquad \text{for all} \quad \xi\in\Omega_k.
\eeq
\er
\noindent
From now on we will often suppress this dependence on the parameter $\rho$, $D(\cdot,\cdot)=D_\rho(\cdot,\cdot)$,  in order not to overload the notation. 
\noindent
For each of the processes we are considering, the orthogonal duality polynomials are given as follows.
\begin{description}
\item[IRW: Charlier polynomials.] The duality polynomials are given by 
\[
d(m,n)=C(m,n)
\]
where $C \(m, \cdot\)$ is the Charlier polynomial of degree $m$ that we characterize by means of its generating function:
\be\label{CharGenFun}
\sum_{m=0}^{\infty} {C(m,n)}\cdot  \frac{t^m}{m!} =  e^{-t} \( \frac{\rho - t}{ \rho} \)^n.
\ee
\item[SEP$(\alpha)$: Krawtchouk polynomials.] Strictly speaking these polynomials do not satisfy a self-duality relation. However, under a proper normalization we can find a duality function in terms of them. The duality polynomials are  given by
\[
d(m,n)=\frac{m! (\alpha- m)!}{\alpha!}\cdot K(m,n)% =\frac {m!} {\pi(m)} \cdot K(m,n)
\]
where  $K(m,\cdot)$ is the Krawtchouk polynomial of degree $m$ whose  generating function is
\be\label{KraGenFun}
\sum_{m=0}^{\infty} K(m,n) \cdot t^m =  \( 1+ t \)^{\alpha} \( \frac{ 1 - \frac{\alpha t}{\rho}}{1+t} \)^n.
\ee
\item[SIP$(\alpha)$: Meixner polynomials.] In this case   the polynomials  satisfying the self-duality relation are given by the following normalization of the Meixner polynomials
\be\label{SIPdualMeix}
d(m,n)=\frac{\Gamma(\alpha)}{\Gamma(\alpha+m)}\cdot M(m,n)  %= \frac 1 {\pi(m)} \cdot M(m,n)
\ee
where $M(m,\cdot)$ is the Meixner polynomial of degree $m$ with generating function
\be\label{MeixGenFun}
\sum_{m=0}^{\infty}  {M(m,n)}\cdot  \frac{t^m}{m!} =  \( 1- t \)^{-\alpha} \( \frac{ 1 - \frac{(\alpha+\rho) t}{ \rho}}{1-t} \)^n.
\ee
\end{description}
We refer the reader to \cite{koekoek1996askey} and \cite{chihara2011introduction} for more details on these polynomials and their generating functions.
%In particular around page 163 in chihara2011
\vskip.1cm

\section{Fluctuation fields}\label{FF}

The density fluctuation field $\caX$ is the stochastic  object usually defined to study small deviations from the hydrodynamic limit. This field corresponds to a central limit type of rescaling of the density field, i.e.
\be\label{simpleflucfield}
\caY_t^{(n)}(\phi,\eta) := {n^{-d/2}}\sum_{ x \in \Zd} \phi(x/n) (\eta_x(n^2t) -\rho).
\ee
\noindent
 Fields of this type have been intensively studied in the literature. For different models,  the sequence $\caY^{(n)}_t$ is proven to converge to  a limiting field $\caY_t$ that  is identified as the distribution-valued random variable satisfying the following martingale problem:
for any $\phi \in S(\R^{d})$ the process 
\be\label{martingX1}
M_t(\phi)= \caY_t(\phi)-\caY_0(\phi) -\frac{\chi \alpha}{2}  \int_0^t  \caY_s(\Delta \phi) ds
\ee
is a square integrable continuous martingale of quadratic variation given by the expression
\be\label{QuadvarmartingX1}
(M_t(\phi))^2 - \chi \rho (\alpha + \sigma \rho) \norm{\nabla \phi(x)}^2 \cdot t
\ee
being also a martingale.\\

\noindent
Formally speaking, \eqref{martingX1} and \eqref{QuadvarmartingX1} imply that the limiting field $\caY_t$ satisfies the Ornstein-Uhlenbeck equation
\be\label{OUSPDE0}
d \caY_t = \tfrac{\chi \alpha}{2} \, \Delta \caY_t \;dt +  \sqrt{ \chi \rho (\alpha + \sigma \rho)} \, \nabla d \caW_t
\ee
where $\caW_t(x)$ is a space-time white noise with covariance
\be
cov [ \caW_t(x), \caW_s(y) ]= \min(t,s) \delta(x-y) \nn 
\ee
%$\cs$ given by
%\be\label{cs}
 %\cs=
%\left\{
%\begin{array}{lll}
%1/ \rho & \text{if} &\sigma=0\\
%{\alpha}/{\rho} & \text{if}& \sigma=1\\
%(\alpha+\rho)/{\rho} & \text{if} & \sigma=-1
%\end{array}
%\right.
%\ee
We refer the reader to \cite{de1984survey} for a precise statement on the convergence for the case of the exclusion process, corresponding, in our setting, to the case  $\alpha=1$ and $\sigma=-1$.\\

\noindent
The density field \eqref{simpleflucfield} can be written, in our context, in terms of our orthogonal polynomial dualities $D_\rho(\xi,\eta)$ by choosing $\xi\in\Omega_1$. Indeed, in all models considered  we have that  there exists a constant $c_{\sigma,\alpha,\rho}$ such that
\beq\label{one}
D_\rho(\delta_x,\eta)= c_{\sigma,\alpha,\rho} \( \eta_x-\rho \)
\eeq
where
\be\label{cs}
 c_{\sigma,\alpha,\rho} =
\left\{
\begin{array}{lll}
1/ \rho & \text{if} &\sigma=0\\
{\alpha}/{\rho} & \text{if}& \sigma=1\\
(\alpha+\rho)/{\rho} & \text{if} & \sigma=-1
\end{array}
\right..
\ee 
Later on, in order not to overload notation we will suppress the dependence on $\rho$ and $\alpha$ and just write $\cs$. From \eqref{one}  we  observe that the field \eqref{simpleflucfield} can be re-written (modulo a multiplicative constant)  as
\be
\caY_t^{(n,1)}(\phi) = {n^{-d/2}} \sum_{ x \in \Zd} \phi\(\tfrac{x}{n}\) D_\rho(\delta_x , \eta(n^2 t))
\ee
where the superindex ${(n,1)}$ suggests that, in some sense, this is the first-order density field. Using \eqref{one} and \eqref{OUSPDE0} the formal limiting SPDE for $\caY_t$ is
\be\label{OUSPDE}
d \caY_t = \frac{\chi \alpha}{2} \, \Delta \caY_t dt +  \cs \sqrt{ \chi \rho (\alpha + \sigma \rho)} \, \nabla d \caW_t
\ee

\noindent
The observation that the field \eqref{simpleflucfield} can be expressed in terms of duality polynomials opens the possibility of defining higher order fields and study their scaling limits. For $k\in \N$, $k\ge 1$ we define the $k$-th order  field as
\beq\label{field}
\caY^{(n,k)}(\phi^{(k)},\eta) := \caX^{(n,k)}(\Phi,\eta) &:=& n^{-kd/2} \sum_{\xi\in\Omega_k}\(\prod_{x\in \Zd} \phi\(\tfrac x n\)^{\xi_x} \)\La(\xi) \cdot D_\rho(\xi,\eta)\\
&=&  n^{-kd/2} \sum_{\xi\in\Omega_k}\(\prod_{x\in \Zd} \phi\(\tfrac x n\)^{\xi_x} \cdot \la(\xi_x) \cdot d_\rho(\xi_x,\eta_x)\)\nn
\eeq
where $\phi\in S(\R^d)$ is a test function, $\La$ is as in \eqref{muk}, and
\beq\label{productphitensor0}
\phi^{(k)}:=\mathop{\bigotimes}_{i=1}^{k} \phi
\eeq
\beq\label{Phi}
\Phi(\xi)=\prod_{x\in \Zd} \phi(x)^{\xi_x}, \qquad \qquad  \Phi_n(\xi)=\prod_{x\in \Zd} \phi\(\tfrac x n\)^{\xi_x}
\eeq
Notice that the only difference between $\caY^{(n,k)}(\phi^{(k)},\eta)$ and $\caX^{(n,k)}(\Phi,\eta)$ is that the latter works on test functions over configuration space, i.e., $\Phi \in S(\Omega_k)$, while the former works on test functions $\phi^{(k)} \in S(\R^{kd})$.
Then, using the notation
\beq
\Di_\rho(\xi,\eta):= \La(\xi) \cdot D_\rho(\xi,\eta), \qquad  \di_\rho(m,n) = \la(m)\cdot d_\rho(m,n)\label{DL}
\eeq
\beq\label{prod}
\Di_\rho(\xi,\eta)=\prod_{i\in \Zd} \di_\rho(\xi_i,\eta_i)
\eeq
we can rewrite  the $k$-th order field \eqref{field} as
\be
\caX^{(n,k)}(\Phi,\eta) := n^{-kd/2} \sum_{\xi\in\Omega_k} \Phi_n(\xi)\cdot  \Di_\rho(\xi,\eta)
\ee
and  define:
\be\label{defkorder3}
\caX_t^{(n,k)}(\Phi) := \caX^{(n,k)}(\Phi,\eta(n^2 t)).
\ee
The choice of multiplying the duality function by the measure $\Lambda(\cdot)$ in \eqref{DL} is dictated simply by computational convenience that, even if obscure at the moment, will be made clearer in the course of the paper.

\section{The coordinate process}

Thinking of $k \in \N$ as the number of particles in our process, we want to introduce a family of permutation-invariant coordinate processes $\{ X^{(k)}(t): t \geq 0 \}$ compatible with the finite configuration processes  $\{\xi(t):t\geq 0\}$ on $\Omega_k$. Here the coordinate process is a Markov process on $\Z^{dk}$ with 
\beq
X^{(k)}(t)= (X_1(t), \ldots, X_k(t)), \qquad X_i(t)\in \Zd, \quad \forall i=1, \ldots, k
\eeq
 $X_i(t)$ being the position of the $i$-th particle at time $t\ge 0$. For a further explanation of the notion of compatibility we refer the reader to \cite{carinci2019consistent}.\\
\noindent
Denote by $\mathbf x \in \mathbb Z^{kd}$ the coordinate vector $\mathbf x:=(x_1,\ldots,x_k)$, with $x_i\in {\mathbb Z}^d$, for $i=1,\ldots, k$. The coordinate process $\{ X^{(k)}(t): t \geq 0 \}$ is defined by means of its infinitesimal generator:
\be\label{gencoord}
L^{(k)} f(\mathbf x) = \sum_{i=1}^{k} \sum_{r \in \caR} p(r) \Bigg( \alpha + \sigma \sum_{\substack{j=1\\ j \neq i}}^k \mathbf 1_{x_j=x_i +r} \Bigg) \( f(\mathbf x^{i,i+r}) -f(\mathbf x) \)
\ee
where $\mathbf x^{i,i+r}$ denotes $\mathbf x$ after moving the particle in position $x_i$ to position $x_i+r \in \Zd$.
\noindent
Notice that for $\mathbf x \in \mathbb Z^{kd}$ the compatible configuration $\xi(\mathbf x) \in \Omega_k$ is given by
\be
\xi(\bix)=\(\xi_i(\bix), i\in \Zd\)\qquad \text{with} \qquad \xi_i(\mathbf x)=\sum_{j=1}^k\mathbf 1_{x_j=i}.
\ee

\subsection{Product $\sigma$-finite reversible measures}

It is possible to verify, by means of detailed balance, that the  coordinate-process $\{X^{(k)}(t) : t\ge 0\}$ admits a reversible  $\sigma$-finite measure 
that is given by
\beq\label{pf}
 \pii(\bix)=   \frac {\La(\xi(\bix))} {N(\xi(\bix))}=  \prod_{i\in \Zd}\xi_i(\bix)!\cdot \la(\xi_i(\bix)) \qquad \text{for} \quad \bix\in\Z^{kd}
\eeq
 with
\beq
N(\xi):=|\{\bix \in \Z^{kd}: \; \xi(\bix)=\xi\}|=\frac{k!}{\prod_{i\in \Zd}\xi_i!}
\eeq
then we can rewrite $\Pi$ in 
 the product form:
\be\label{defpirho}
\pii(\mathbf x) =  \mathop{\prod}_{i \in \Zd } \pi(\xi_i(\mathbf x)), \qquad \quad  \bix =(x_1,\ldots, x_k)\in \Z^{kd}
\ee
with $\pi$ given as follows:
\be\label{defpirhoIRWmarg}
\pi(m)= m! \cdot \la(m)=
\left\{
\begin{array}{llll}
 1, &m\in \N& \text{for }\sigma=0& \text{IRW}\\
 \\
 \frac{\alpha!}{(\alpha-m)!}, & m \in \{0,\ldots, \alpha\} & \text{for }\sigma=-1& \text{SEP}(\alpha)\\
 &\\
\frac{\Gamma(\alpha + m)}{\Gamma(\alpha)},& m\in \N &\text{for }\sigma=1 & \text{SIP}(\alpha)
\end{array}
\right..\ee
\noindent
%Notice that $\tfrac 1 {k!}\cdot\Pi$ is a probability measure on $\Omega_k$, indeed, from \eqref{pf} we have
%\be
%\frac 1 {k!}\sum_{\bix \in \Omega_k} \Pi(\bix)=1
%\ee
%\noindent
Given the measures $\pii$, we now consider the spaces of permutation-invariant functions:
\be
\hat{L}^2(\mathbb Z^{kd}, \pii):= \Big\{ f \in L^2(\mathbb Z^{kd}, \pii) : f(\bix) = f(\bix_{\sigma}), \: \:\forall \sigma \in \caP(k) \Big\}
\ee
with $\caP(k)$ denoting the set of all possible permutations of the set $\{ 1,2, 3, \dots k\}$. We  endowed the space $\hat{L}^2(\mathbb Z^{kd}, \pii)$ with the inner product given by:
\be
\langle f, g \rangle_{\pii} = \sum_{\bix \in \mathbb Z^{kd}} f(\bix) g(\bix) \pii(\bix).
\ee

\br
Notice that any function $f \in \hat{L}^2(\mathbb Z^{kd}, \pii)$ can be interpreted also as a function on the configuration space. In this work we will extensively use this fact  by changing between interpretations sometimes from one line to another in the same derivation.
\er

\br\label{adjointlock}
As a consequence of reversibility of the measures $\pii$, we can infer that the $k$-particles generator $L^{(k)}$ is self-adjoint with respect to the inner product $\langle\cdot,\cdot\rangle_{\pii}$, i.e.
\be
\langle f, L^{(k)} g \rangle_{\pii} = \langle L^{(k)}f, g \rangle_{\pii}
\ee
for all $f, g \in \hat{L}^2(\mathbb Z^{kd}, \pii)$.
\er

\subsection{The fluctuation fields in  coordinate notation}
It is possible to rewrite the fluctuation field \eqref{field} in the coordinate variables.
Notice that in this context the test function $\Phi$ defined in \eqref{Phi} becomes a tensor function:
\beq
\Phi(\xi(\bix))= \prod_{i=1}^n \phi(x_i)
\eeq
i.e. it is the  homogeneous $k$-tensor test function $\phi^{(k)}\in S(\R^{kd})$ of the form
\beq\label{productphitensor}
\Phi \circ \xi=\phi^{(k)} =\mathop{\bigotimes}_{i=1}^{k} \phi
\eeq
then, after a change of variable in the sum we can rewrite  the $k$-th field  as follows  
\be\label{defkorder}
\caY^{(n,k)}(\phi^{(k)},\eta) =  \caX^{(n,k)}(\Phi,\eta) =n^{-kd/2} \sum_{\bix \in \Z^{kd}} \phi^{(k)}\(\tfrac{\bix}{n}\) \cdot  \pii(\bix)\cdot D(\xi(\bix), \eta).
\ee
Notice that we can also let the field $\caY$ act on a general $f \in S(\R^{kd})$ as expected, i.e.,
\be\label{defkorderNONSYM}
\caY^{(n,k)}(f,\eta) = n^{-kd/2} \sum_{\bix \in \Z^{kd}} f\(\tfrac{\bix}{n}\) \cdot  \pii(\bix)\cdot D(\xi(\bix), \eta).
\ee
\br
Because we deal with unlabeled particle systems it is natural to define the higher order fluctuation fields acting on symmetric test functions $\Phi$ i.e. on elements of the Schwartz space  $S(\R^{kd})$ that are permutation-invariant: $\Phi(x_{\sigma(1)},\ldots,x_{\sigma(k)})=\Phi(x_1,\ldots, x_k)$ for all $\sigma\in \mathcal{P}(k)$, the set of permutations of $\{1,\ldots, k\}$. By polarization  and the fact that the orthogonal polynomials are functions of the unlabeled particle configurations, it is therefore sufficient to restrict our analysis to the set of symmetric tensor products of $S(\R^d)$. Indeed the linear combinations of those generate the symmetric test functions, see \cite{thomas2014polarization}.
\er

\section{Main result}

\subsection{Heuristics: macroscopic dynamics}
The goal of this section is to provide some intuitions on the type of limiting field that we should expect for fields of order greater than one. We will start by considering the cases $k=1,2$ and,  inspired by the results obtained in \cite{gonccalves2019quadratic},
we will propose a heuristic  interpretation of the two  SPDEs obtained as scaling limits and their relation. Based on this interpretation we  will conjecture  a possible generalization to the $k$th-order  case. In Section \ref{MT} we will give the rigorous result confirming the validity of the conjecture.
\vskip.3cm
\noindent
Here we will informally use the notation $\caX^{(k)}_t$ and $\caY^{(k)}_t$ for the distributional limits of $\caX^{(n,k)}$ and $\caY^{(n,k)}$ respectively. 
\vskip.3cm
%\noindent
Recall that from \eqref{OUSPDE} we know  that formally the distribution valued first order field $\caY^{(1)}_t(x)$ is a solution to the Ornstein-Uhlenbeck equation
\be\label{OUSPDEheuristics}
d \caY^{(1)}_t(x) = \frac{\chi \alpha}{2} \, \Delta \caY^{(1)}_t(x) \; dt + \cs\sqrt{ \chi \rho (\alpha + \sigma \rho)} \, \nabla d\caW_t(x)
\ee
where for $ x \in \Rd$, $\caW_t(x)$ is a space-time white noise.
%with covariance
%\be
%cov [ \caW_t(x), \caW_s(y) ]= \min(t,s) \delta(x-y) \nn 
%\ee
\noindent
Additionally, from the martingale problem given in \cite{gonccalves2019quadratic}, we can deduce that the distribution-valued second-order field $\caY^{(2)}_t(x,y)$ is a solution to the SPDE
\beq\label{2ndSPDEheuristics}
d \caY^{(2)}_t(x,y) &=& \frac{\chi \alpha}{2} \, \Delta^{(2)} \caY^{(2)}_t(x,y) dt + \cs\sqrt{ \chi \rho (\alpha + \sigma \rho)} \, \caY^{(1)}_t(x) \nabla d\caW_t(y) \nn \\
&+&  \cs\sqrt{ \chi \rho (\alpha + \sigma \rho)} \, \caY^{(1)}_t(y) \nabla d \caW_t(x)
\eeq
where $\caW_t(x)$ is the white noise in \eqref{OUSPDEheuristics} and  $\Delta^{(2)}$ is the $\R^{2d}$ dimensional Laplacian, which is the sum of the Laplacian in the $x$ variable plus the Laplacian in the $y$ variable.\\

The key idea to extrapolate these relations to higher orders is to interpret the non-linearity on the RHS of \eqref{2ndSPDEheuristics} as a fields product, that we denote by $\diamond$, that satisfies the Leibniz rule of differentiation. This interpretation suggests that the second-order field $\caY^{(2)}_t(x,y)$ is, in turn, a second power of the first-order field $\caY^{(1)}_t(x)$. More precisely conjecturing
\be
\caY^{(2)}_t(x,y) = \caY^{(1)}_t(x) \diamond \caY^{(1)}_t(y), \nn
\ee
since the product $\diamond$ follows the Leibniz rule we would  have that
\beq\label{wick2ndSPDEheuristics}
d \caY^{(2)}_t(x,y) &=& d \(  \caY^{(1)}_t(x) \diamond \caY^{(1)}_t(y) \) \nn \\
&=& d \caY^{(1)}_t(x) \diamond \caY^{(1)}_t(y) + \caY^{(1)}_t(x) \diamond d \caY^{(1)}_t(y) \nn \\
&=& \(\frac{\chi \alpha}{2} \, \Delta \caY^{(1)}_t(x) dt +\cs\sqrt{ \chi \rho (\alpha + \sigma \rho)} \, \nabla d\caW_t(x)   \)  \diamond \caY^{(1)}_t(y)  \nn \\
&+& \caY^{(1)}_t(x) \diamond \( \frac{\chi \alpha}{2} \, \Delta \caY^{(1)}_t(y) dt + \cs\sqrt{ \chi \rho (\alpha + \sigma \rho)} \, \nabla d\caW_t(y)    \) \nn \\
&=& \frac{\chi \alpha}{2} \, \Delta^{(2)} \caY^{(2)}_t(x,y) dt + \cs\sqrt{ \chi \rho (\alpha + \sigma \rho)} \, \caY^{(1)}_t(x) \diamond \nabla d\caW_t(y) \nn \\
&+&  \cs\sqrt{ \chi \rho (\alpha + \sigma \rho)} \, \caY^{(1)}_t(y) \diamond \nabla d\caW_t(x)
\eeq
which indeed agress with \eqref{2ndSPDEheuristics}.\\

\noindent
After the discussion above, it seems natural to expect that the $k$th-order field is a $k$th $\diamond$-power of the first order one. More precisely
we conjecture that a relation of the type
\be
\caY^{(k)}_t(x_1,x_2,\dots,x_k) = \caY^{(1)}_t(x_1) \diamond \caY^{(1)}_t(x_2) \diamond \dots \diamond \caY^{(1)}_t(x_k).  \nn
\ee
is satisfied. If this holds true,  computations analogous  to \eqref{wick2ndSPDEheuristics}  would imply the formal SPDE
\beq\label{wickkthSPDEheuristics}
d \caY^{(k)}_t(\mathbf x) &=& \frac{\chi \alpha}{2} \, \Delta^{(k)} \caY^{(k)}_t(\mathbf x) dt + \cs \sqrt{ \chi \rho (\alpha + \sigma \rho)} \, \sum_{j=1}^k \caY^{(k-1)}_t(\mathbf x^{-j}) \diamond \nabla d\caW_t(x_j) 
\eeq
where $\Delta^{(k)}$ is the $\R^{kd}$-dimensional Laplacian, defined as the sum of the Laplacians at each coordinate and $\mathbf x^{-j}$ is the $(k-1)d$-dimensional vector  obtained from $\mathbf x$ by removing its coordinate $x_j$.\\
\vskip.3cm
\noindent
The equation that we thus informally obtained  will be later on justify in the main result of the next section.

\subsection{Main theorem}\label{MT}
Let us spend one paragraph to introduce the probability notions which are relevant for our main result. As we already mentioned, the $k$th-order fluctuation field can be considered as taking values in $S\myprime(\R^k)$, the space of tempered distributions which is dual to $S(\R^k)$. Our original process $\eta_{n^2t}$ has state space $\Omega^{(n)}$  corresponding to the rescaled lattice $\tfrac{1}{n}\Z$. We then denote by $\P_n$, respectively $\E_n$, the probability measure, respectively expectation, induced by the measure $\nu_{\rho}$ and the diffusively rescaled process $\eta_{n^2t}$ on  $D([0,T];\Omega^{(n)})$. Finally, we denote by $Q_n^{(k)}$ the probability measure on $D([0,T];S\myprime(\R^k))$ induced by the density fluctuation field $ \caY_t^{(n,k)}$ over $\P_n$. 
\vskip.3cm
\noindent

\vskip.3cm
\noindent
\bt\label{mainkthorder}
The process $\{\caY_{t}^{(n,k)} : t \in [0,T] \}$ converges in distribution, with respect to the J1-topology of $D([0,T];S\myprime(\R^{kd}))$, as $n \to \infty$ to the process $\{\caY_t^{(k)} : t \in [0,T] \}$ being the unique solution of the following recursive martingale problem.
\vskip.3cm
\noindent
{\bf{Recursive martingale problem:}} for any symmetric $\phi^{(k)} \in S(\R^{kd})$ the process 
\be\label{martingXk}
M_t^{(k)}(\phi^{(k)})= \caY_t^{(k)}(\phi^{(k)})-\caY_0^{(k)}(\phi^{(k)}) -\frac{\chi \alpha}{2} \, \int_0^t  \caY_s^{(k)}(\Delta^{(k)}\phi^{(k)}) ds
\ee
is a continuous square integrable martingale of quadratic variation
\be\label{QuadvarmartingXk}
(M_t^{(k)}(\phi^{(k)}))^2 -\cs^2\chi \rho (\alpha + \sigma \rho) \int_0^t  \int_{\R^d} \norm{\nabla \phi(x)}^2 \( \caY_s^{(k-1)} ( \phi^{(k-1)}) \)^2 dx \, ds
\ee
with initial condition $\caY_t^{(1)}$ given by the solution of \eqref{OUSPDEheuristics}.
\et
\br
This recursive martingale problem is the rigorous counter part of the formal SPDE \eqref{wickkthSPDEheuristics} that we heuristically obtained.
\er

\subsection{Strategy of the proof}
\noindent
We will show Theorem \ref{mainkthorder} by using induction on $k$. In the proof we will take advantage of the fact that the base case, $k=1$, is already proved in the literature. On the other hand, the inductive step will be proven by means of an approach based on the natural Dynkin martingales: 
\be\label{Mtdef}
M_t^{(n,k)} (\Phi) = \caX_{t}^{(n,k)}(\Phi) - \caX_0^{(n,k)}(\Phi) - n^2 \int_0^t \loc \caX_{s}^{(n,k)}(\Phi) ds
\ee
and
\be\label{Ntdef}
N_t^{(n,k)} (\Phi) = (M_t^{(n,k)} (\Phi))^2 - n^2 \int_0^t \Gamma \caX_{s}^{(n,k)}(\Phi) ds
\ee
where $\Gamma$ is the so-called carr\'e-du-champ operator given by:
\be\label{CdC}
\Gamma (f) = \loc (f^2) - 2 f \loc (f).
\ee
Notice that the Dynkin martingales can also be expressed in terms of the fields $ \caY_t^{(n,k)}$.\\

\noindent
Roughly our approach consists of the following steps:
\ben
\item we express the integrand term of equation \eqref{Mtdef} in terms of the $k$th-order fluctuation field $\caX^{(n,k)}$ using duality (Section \ref{Drift});
\item we close the  equation \eqref{Ntdef} by expressing the integrand in the RHS in terms of the $(k-1)$th-order fluctuation field $\caX^{(n,k-1)}$ (Section \ref{CdC});
\item we show tightness for the sequence of probability measures $Q_n^{(k)}$ (Section \ref{T});
\item finally we characterize the limiting field by showing uniqueness of the solution of the martingale problem (Sections  \ref{Char}-\ref{Uni}).
\een

\subsection{Inductive argument}
The proof is done by induction over the order of the field $k$.
 The base case  $k=1$, corresponding to the density fluctuation field \eqref{simpleflucfield}, is assumed to be true. Indeed, as  mentioned in Section \ref{FF},  a proof of Theorem \ref{mainkthorder} for exclusion dynamics  and zero-range processes (of which independent random walkers are a particular case) is given in \cite{de1984survey} and \cite{kipnis2013scaling} respectively. By similar arguments the  result can be extended to the case of inclusion process.

\vskip.3cm
\noindent
To implement the inductive argument we formalize the following inductive hypothesis that will be referred to several times in the course of the proof of Theorem \ref{mainkthorder}.
\bind\label{InducHyp}
For any $k_0 \in \{1,2,\dots,k-1 \}$ the sequence $\{\caY_{t}^{(n,k_0)} : t \in [0,T] \}$ converges in distribution, with respect to the J1-topology of $D([0,T];S\myprime(\R^{k_0d}))$, as $n \to \infty$ to the process $\{\caY_t^{(k_0)} : t \in [0,T] \}$ being the unique solution of the following martingale problem.
\vskip.3cm
\noindent
{\bf{Martingale problem:}}
for any symmetric $\phi^{(k_0)} \in S(\R^{k_0d})$ the process 
\be\label{martingXk_0}
M_t^{(k_0)}(\phi^{(k_0)})= \caY_t^{(k_0)}(\phi^{(k_0)})-\caY_0^{(k_0)}(\phi^{(k_0)}) -\frac{\chi \alpha}{2} \, \int_0^t  \caY_s^{(k_0)}(\Delta^{(k_0)}\phi^{(k_0)}) ds
\ee
is a continuous square integrable martingale of quadratic variation
\be\label{QuadvarmartingXk_0}
(M_t^{(k_0)}(\phi^{(k_0)}))^2 -  \cs^2\chi \rho (\alpha + \sigma \rho) \int_0^t  \int_{\R^d} \norm{\nabla \phi(x)}^2 \( \caY_s^{(k_0-1)} ( \phi^{(k_0-1)}) \)^2 dx \, ds.
\ee
\eind
\noindent
%We now follow the steps 1-4 sketched in the previous section to show that Theorem \ref{mainkthorder} is valid for $k$.

{\section{Proof of Theorem \ref{mainkthorder}}}

\subsection{Closing the equation for the drift term: $k \geq 2$}\label{Drift}

In order to close the equation for the drift term  thanks to Remark \ref{adjointlock} we can just proceed as follows
\beq\label{splitgenkparticles}
n^2 \loc \caX^{(n,k)}(\Phi,\eta) &=& n^{-kd/2} \sum_{\xi \in \Omega_k} n^2\Phi_n(\xi)\cdot  [\loc \Di(\xi, \cdot)](\eta) \nn\\
&=& n^{-kd/2} \sum_{\xi \in \Omega_k} n^2\Phi_n(\xi)\cdot \La(\xi) \cdot [\loc D(\xi, \cdot)](\eta)\nn\\
&=& n^{-kd/2} \sum_{\xi \in \Omega_k} n^2\Phi_n(\xi)\cdot \La(\xi) \cdot [\loc^{(k)} D(\cdot, \eta)](\xi)\nn\\
&=& n^{-kd/2} \sum_{\xi \in \Omega_k} n^2 [\loc^{(k)} \Phi_n](\xi)\cdot \La(\xi) \cdot D(\xi, \eta)\nn\\
&=& n^{-kd/2} \sum_{\xi \in \Omega_k} n^2 [\loc^{(k)} \Phi_n](\xi)\cdot  \Di(\xi, \eta).\nn
\eeq
We proceed  evaluating the action of the $k$-particles generator on $\Phi_{n}$. We then have 
\beq\label{kgeq2sigma}
n^2 [\loc^{(k)} \Phi_n](\xi) &=&  \sum_{x \in \Zd} \sum_{ r \in \caR} p(r) \cdot \xi_x ( \alpha + \sigma \xi_{x+r})\cdot n^2(\Phi_n(\xi^{x,x+r})- \Phi_n(\xi))\nn\\
 &=&\sum_{x \in \Zd}\Phi_n(\xi-\delta_x)\sum_{ r \in \caR} p(r) \cdot \xi_x ( \alpha + \sigma \xi_{x+r})\cdot n^2\(\phi\(\tfrac{x+r}n\)- \phi\(\tfrac xn\)\)\nn\\
&=&\alpha\sum_{x \in \Zd} \Phi_n(\xi-\delta_x) \cdot \xi_x \sum_{ r \in \caR} p(r) \cdot n^2\(\phi\(\tfrac{x+r}n\)- \phi\(\tfrac xn\)\)\nn\\
&+&\sigma\sum_{x \in \Zd} \Phi_n(\xi-\delta_x)\sum_{ r \in \caR} p(r) \cdot \xi_x  \xi_{x+r}\cdot n^2\(\phi\(\tfrac{x+r}n\)- \phi\(\tfrac xn\)\).
\eeq
\br
Notice that the contribution coming from the second term in the RHS of \eqref{kgeq2sigma} does not appear in the case $k=1$.
\er
\noindent
First of all we prove that
\be\label{taylorexpankfirst}
n^2 \sum_{r \in \caR} p(r)  \(\phi\(\tfrac{x+r}n\)- \phi\(\tfrac xn\)\) =\tfrac \chi 2\cdot \Delta{\phi}\(\tfrac{x}n\) +  \tfrac 1 n \psi_{n}\(\tfrac {x}n\)
\ee
for a suitable $\psi_{n}\in S(\R)$ such that 
\beq\label{QUO}
\sup_n\frac 1 {n^d} \sum_{x\in\Zd} \psi_{n}\(\tfrac {x}n\)<\infty.
\eeq
To prove this we use the Taylor expansion:
\beq\label{taylor}
\phi\(\tfrac{x+r}n\) -\phi\(\tfrac{x}n\)=\frac 1 n \sum_{j=1}^d r_j  \cdot \frac{\partial \phi}{\partial{x_{j}}}\(\tfrac{x}n\)+\frac 1 {2n^2}\sum_{j,\ell=1}^d r_j r_\ell \frac{\partial^2 \phi}{\partial{x_{j}} \partial x_{\ell}}\(\tfrac{x}n\)+ \ldots
\eeq
and then
\beq
&&n^2 \sum_{r \in \caR} p(r)  \(\phi\(\tfrac{x+r}n\)- \phi\(\tfrac xn\)\) \nn\\
&&=n\sum_{j=1}^d \(\sum_{r \in \caR}r_j  p(r)\)   \cdot \frac{\partial \phi}{\partial{x_{j}}}\(\tfrac{x}n\)+\frac 12 \sum_{j,\ell=1}^d \(\sum_{r \in \caR} r_j r_\ell p(r)\) \frac{\partial^2 \phi}{\partial{x_{j}}\partial x_{\ell}}\(\tfrac{x}n\)+  \ldots \nn
\eeq
for some $\psi_{n}$ satisfying \eqref{QUO}.
From the assumption \eqref{ps}, it follows that:
\beq
\sum_{r_j=-R}^R r_j  p(r)=0
\eeq
thus, from  the fact that $\caR=[-R,R]^d\cap \Zd$ we have
\beq
\sum_{r\in \caR} r_j  p(r)=0\qquad \text{and} \qquad \sum_{r \in \caR} r_j r_\ell p(r)=0 \quad \text{for} \:j \neq \ell
\eeq
as a consequence,
\beq
&&n^2 \sum_{r \in \caR} p(r)  \(\phi\(\tfrac{x+r}n\)- \phi\(\tfrac xn\)\) =\frac 1 2 \sum_{\ell=1}^d \(\sum_{r \in \caR}  r^2_\ell p(r)\) \frac{\partial^2 \phi}{\partial x^2_{\ell}}\(\tfrac{x}n\)+   \tfrac 1 n \psi_{n}\(\tfrac {x}n\)
\nn\\
&&=\tfrac\chi 2 \cdot \sum_{\ell=1}^d \frac{\partial^2 \phi}{\partial x^2_{\ell}}\(\tfrac{x}n\)+   \tfrac 1 n \psi_{n}\(\tfrac {x}n\)
\nn
\eeq
from which it follows \eqref{taylorexpankfirst}. 
\vskip.2cm
\noindent
Now we have
\beq
n^2 [\loc^{(k)} \Phi_n](\xi)  &=&\alpha\sum_{x \in \Zd} \Phi_n(\xi-\delta_x) \cdot \xi_x\cdot \(\tfrac \chi 2\cdot \Delta{\phi}\(\tfrac{x}n\) +   \tfrac 1 n \psi_{n}\(\tfrac {x}n\)\)+E_n(\phi,\xi)
\nn
\eeq
with
\beq\label{EEo}
E_n(\phi,\xi)  &:=& \sigma\sum_{x \in \Zd}\Phi_n(\xi-\delta_x)\sum_{ r \in \caR} p(r) \cdot \xi_x  \xi_{x+r}\cdot n^2 \(\phi\(\tfrac{x+r}n\)- \phi\(\tfrac xn\)\)
\eeq
then we have
\beq\label{ciao}
&&\loc \caX^{(n,k)}(\Phi,\eta)- \frac{1}{n^{kd/2}} \sum_{\xi\in \Omega_k} E_n(\phi,\xi)   \cdot   \Di(\xi, \eta)  \nn\\
&&=  \frac{\alpha}{n^{kd/2}} \sum_{\xi\in \Omega_k}\Di(\xi,\eta)\sum_{x \in \Zd} \Phi_n(\xi-\delta_x) \cdot \xi_x\cdot \(\tfrac \chi 2\cdot \Delta{\phi}\(\tfrac{x}n\)  +   \tfrac 1 n \psi_{n}\(\tfrac {x}n\)\).\nn
\eeq
It is  now convenient to pass to the coordinate notation to treat sums of the type:
\beq
&& \sum_{\xi\in \Omega_k}\Di(\xi,\eta) \sum_{x \in \Zd} \Phi_n(\xi-\delta_x) \cdot \xi_x\cdot \psi(\tfrac x n)\nn
\eeq
for some  $\psi\in S(\R^d)$. First of all we notice that summing over $\xi\in\Omega_k$ is the same as summing over $\bix\in\Z^{kd}$:
\beq
&& \sum_{\xi\in \Omega_k}\Di(\xi,\eta) \sum_{x \in \Zd} \Phi_n(\xi-\delta_x) \cdot \xi_x\cdot \psi(\tfrac x n)\nn\\
&=& \sum_{\bix\in \Z^{kd}}\frac 1{N(\xi(\bix))}\cdot \Di(\xi(\bix),\eta) \sum_{i=1}^k \Phi_n(\xi(\bix)-\delta_{x_i}) \cdot \psi(\tfrac {x_i} n)\nn\\
&=& \sum_{\bix\in \Z^{kd}}\frac {\La(\xi(\bix))}{N(\xi(\bix))}\cdot D(\xi(\bix),\eta) \sum_{i=1}^k \psi(\tfrac {x_i} n) \prod_{\substack{\ell=1\\ \ell\neq i}}^k \phi(\tfrac{x_\ell}n)\nn\\
&=& k \sum_{\bix\in \Z^{kd}}\Pi(\bix)\cdot D(\xi(\bix),\eta)  \prod_{\ell=1}^{k-1} \phi(\tfrac{x_\ell}n)\cdot \psi(\tfrac {x_k} n) \nn\\
&=& k  n^{kd/2}\;\caY^{(n,k)}(\phi^{(k-1)}\otimes \psi,\eta)\nn
\eeq
where the last identity follows using the expression of the field acting on more general (i.e., non-symmetric)  test functions  \eqref{defkorderNONSYM}.
Then, substituting in \eqref{ciao} we get
\beq\label{ciao1}
&&\loc \caX^{(n,k)}(\Phi,\eta)- \frac{1}{n^{kd/2}} \sum_{\xi\in \Omega_k} E_n(\phi,\xi)   \cdot   \Di(\xi, \eta)  \nn\\
&=& \alpha k \caY^{(n,k)}\(\phi^{(k-1)}\otimes (\tfrac \chi 2 \Delta \phi+\tfrac 1 n \psi_{n}),\eta\)\nn
\nn
\eeq
where we used the fact that $\phi$ is uniformly bounded on $\Z$. 
From this we can see that it is possible  to close the equation for the second order fluctuation field, modulo an error term that we define as follows
\beq\label{E}
\caE^{(n,k)}(\phi,\eta):=\loc \caX^{(n,k)}(\Phi,\eta)-\alpha k \cdot\tfrac {\chi}2 \cdot \caY^{(n,k)}(\phi^{(k-1)}\otimes \Delta\phi,\eta)
\eeq
\vskip.2cm
\noindent
then we have
\beq\label{EE}
\caE^{(n,k)}(\phi,\eta)= \caE_1^{(n,k)}(\phi,\eta)+\caE_2^{(n,k)}(\phi,\eta)
\eeq
with
\beq
\caE_1^{(n,k)}(\phi,\eta):= \frac{\alpha k}n \;\caY^{(n,k)}\(\phi^{(k-1)}\otimes  \psi_{n},\eta\)\nn
\eeq
and
\beq
\caE_2^{(n,k)}(\phi,\eta):=\frac{1}{n^{kd/2}} \sum_{\xi\in \Omega_k} E_n(\phi,\xi)   \Di(\xi, \eta) 
\eeq
that has to be estimated. Analogously to the previous computation we have
\beq
E_n(\phi,\xi(\bix)) &=&  \sigma n^2 \sum_{i=1}^k \Big(\prod_{\substack{\ell=1\\\ell\neq i}}^k \phi(\tfrac {x_\ell}n) \Big)\cdot
\sum_{r\in\caR} p(r)\(\sum_{j=1}^k \mathbf 1_{x_j=x_i+r}\) \(\phi(\tfrac{x_i+r}n)-\phi(\tfrac{x_i}n)\)\nn\\
&=& \sigma n^2 \sum_{i=1}^k \Big(\prod_{\substack{\ell=1\\\ell\neq i}}^k \phi(\tfrac {x_\ell}n) \Big)\cdot
\sum_{j=1}^k  p(x_j-x_i)\(\phi(\tfrac{x_j}n)-\phi(\tfrac{x_i}n)\)\nn\\
&=& \sigma n^2\sum_{i,j=1}^k \Big(\prod_{\substack{\ell=1\\\ell\neq i,j}}^k \phi(\tfrac {x_\ell}n) \Big)\cdot
 p(x_j-x_i)\phi(\tfrac{x_j}n)\(\phi(\tfrac{x_j}n)-\phi(\tfrac{x_i}n)\)\nn\\
 &=& \sigma  \sum_{\substack{\{i,j\} \\ 1 \leq i, j \leq k}} \Big( \prod_{\substack{\ell=1\\\ell\neq i,j}}^k \phi(\tfrac {x_\ell}n) \Big)\cdot
 p(x_j-x_i)n^2\(\phi(\tfrac{x_j}n)-\phi(\tfrac{x_i}n)\)^2\nn
  %&=& \sigma \sum_{\substack{\{i,j\} \\ 1 \leq i, j \leq k}} \Big( \prod_{\substack{\ell=1\\\ell\neq i,j}}^k \phi(\tfrac {x_\ell}n) \Big)\cdot p(x_j-x_i)\cdot  \Psi_n(\tfrac{x_i}n,\tfrac{x_j}n)\nn
%&&= \sigma n^2 \sum_{i=1}^{k}  \phi_n^{(k-1)}({\mathbf x^{-i}})  \sum_{r \in \caR} p(r)   \sum_{\substack{j=1\\ j \neq i}}^k \mathbf 1_{x_i +r =x_j}    \(\phi\(\tfrac{x_i+r}n\)- \phi\(\tfrac {x_i}n\)\)  \nn \\
%&=& \sigma n^2 \sum_{i=1}^{k}\sum_{j=1}^{k}  \phi_n^{(k-1)}({\mathbf x^{-i}})  p(x_j-x_i)   \(\phi\(\tfrac{x_j}n\)- \phi\(\tfrac {x_i}n\)\) \nn \\
%&=& \sigma n^2 \sum_{i=1}^{k}\sum_{j=1}^{k}  \phi_n^{(k-2)}({\mathbf x^{-i-j}})  \phi_n({x_j}) p(x_j-x_i)     \(\phi\(\tfrac{x_j}n\)- \phi\(\tfrac {x_i}n\)\)  \nn \\
%&=& \sigma n^2 \sum_{\substack{(i,j) \\ 1 \leq i, j \leq k}}  \phi_n^{(k-2)}({\mathbf x^{-i-j}}) p(x_j-x_i)      \(\phi\(\tfrac{x_j}n\)- \phi\(\tfrac {x_i}n\)\)^2 \nn \\
%&=& \sigma \sum_{\substack{(i,j) \\ 1 \leq i, j \leq k}}  \phi_n^{(k-2)}({\mathbf x^{-i-j}}) p(x_j-x_i) \(x_j-x_i\)^2  \(  \frac{\partial  \phi}{\partial x_i} (\tfrac{x_i}n) \)^2 \nn \\
\eeq
where in the last step we used the symmetry of $p(\cdot)$. Then
\beq\label{kparticlessecondterm}
&&\caE_2^{(n,k)}(\phi,\eta)= \nn\\
%&= &\sum_{\bix\in \Z^{kd}}  \frac 1 {N(\xi(\bix))} \cdot \Di(\xi(\bix), \eta)  \cdot E_n(\phi,\xi(\bix))   \nn\\
&&= \frac{1}{n^{kd/2}}   \sum_{\bix\in \Z^{kd}} \Pi(\bix)\cdot D(\xi(\bix), \eta)  \cdot E_n(\phi,\xi(\bix))   \nn\\
&&=\frac{\sigma}{n^{kd/2}}    \sum_{\bix\in \Z^{kd}} \Pi(\bix)\cdot D(\xi(\bix), \eta)  \cdot \sum_{\substack{\{i,j\} \\ 1 \leq i, j \leq k}} \Big( \prod_{\substack{\ell=1\\\ell\neq i,j}}^k \phi(\tfrac {x_\ell}n) \Big)\cdot
 p(x_j-x_i)\cdot n^2\(\phi(\tfrac{x_j}n)-\phi(\tfrac{x_i}n)\)^2\nn\\  \nn\\
 &&=\frac{k(k-1)\sigma}{2n^{kd/2}}   \sum_{\bix\in \Z^{kd}} \Pi(\bix)\cdot D(\xi(\bix), \eta)  \cdot \Big(\prod_{\ell=1}^{k-2}\phi(\tfrac {x_\ell}n) \Big)\cdot
 p(x_{k}-x_{k-1}) \cdot n^2\(\phi(\tfrac{x_{k-1}}n)-\phi(\tfrac{x_k}n)\)^2.\nn
%&&\col{=\frac{k\sigma}{n^{kd/2}}   \sum_{\bix\in \Z^{(k-1)d}}\Big(\prod_{\ell=1}^{k-2}\phi(\tfrac {x_\ell}n) \Big)\cdot  \sum_{r\in\caR}\Pi((\bix,x_{k-1}+r))\cdot D(\xi((\bix,x_{k-1}+r)), \eta)  \cdot \Psi_n(\tfrac{x_{k-1}}n,\tfrac{x_{k-1}+r}n) \cdot p(r)} \nn
%  &&\col{=\frac {k\sigma}{n^{d/2}} \; \caY^{(n,k-1)}(\phi^{(k-2)}\otimes \widetilde \psi_n,\eta) }
%\eeq
%\col{this would be true if one can prove that
% \beq
% &&  \sum_{r\in\caR}\Pi((\bix,x_{k-1}+r))\cdot D(\xi((\bix,x_{k-1}+r)), \eta)  \cdot \Psi_n(\tfrac{x_{k-1}}n,\tfrac{x_{k-1}+r}n) \cdot
% p(r) \nn\\
%  &&\simeq \Pi((\bix,x_{k-1}))\cdot D(\xi((\bix,x_{k-1})), \eta)  \cdot \widetilde \psi_n(\tfrac{x_{k-1}}n) \cdot
\eeq
%??????????????????}
%meaning here, by $O(n^{-1})$, a function of $\bix$ that is uniformly bounded by a constant times $n^{-1}$ for all $\bix\in \Z^2$. 
Hence we have
\beq\label{EEE}
&&\caE^{(n,k)}(\phi,\eta)
%&= &\sum_{\bix\in \Z^{kd}}  \frac 1 {N(\xi(\bix))} \cdot \Di(\xi(\bix), \eta)  \cdot E_n(\phi,\xi(\bix))   \nn\\
= \frac{k}{n^{kd/2}}   \sum_{\bix\in \Z^{kd}} \Pi(\bix)\cdot D(\xi(\bix), \eta)  \cdot \Psi_n(\bix)%&&\col{=\frac{k\sigma}{n^{kd/2}}   \sum_{\bix\in \Z^{(k-1)d}}\Big(\prod_{\ell=1}^{k-2}\phi(\tfrac {x_\ell}n) \Big)\cdot  \sum_{r\in\caR}\Pi((\bix,x_{k-1}+r))\cdot D(\xi((\bix,x_{k-1}+r)), \eta)  \cdot \Psi_n(\tfrac{x_{k-1}}n,\tfrac{x_{k-1}+r}n) \cdot p(r)} \nn
%  &&\col{=\frac {k\sigma}{n^{d/2}} \; \caY^{(n,k-1)}(\phi^{(k-2)}\otimes \widetilde \psi_n,\eta) }
%\eeq
%\col{this would be true if one can prove that
% \beq
% &&  \sum_{r\in\caR}\Pi((\bix,x_{k-1}+r))\cdot D(\xi((\bix,x_{k-1}+r)), \eta)  \cdot \Psi_n(\tfrac{x_{k-1}}n,\tfrac{x_{k-1}+r}n) \cdot
% p(r) \nn\\
%  &&\simeq \Pi((\bix,x_{k-1}))\cdot D(\xi((\bix,x_{k-1})), \eta)  \cdot \widetilde \psi_n(\tfrac{x_{k-1}}n) \cdot
\eeq
with
\beq\label{PSI}
\Psi_n(\bix)&:=& \phi^{(k-2)}(x_1,\ldots,x_{k-2})\otimes \left(\tfrac{\alpha } n\;  \phi(x_{k-1})\cdot\psi_n(\tfrac{x_k}n) \right. \nn \\
&+& \left. \tfrac{\sigma(k-1)}2 p(x_{k}-x_{k-1})  n^2\(\phi(\tfrac{x_{k-1}}n)-\phi(\tfrac{x_k}n)\)^2\right).
\eeq
\noindent
It remains to  show that the $L^{2}(\P_n)$ norm of $\caE^{(n,k)}(\phi, \eta(n^2 t))$ vanishes in the limit as $n$ goes to infinity. This is done in the following lemma:

%\subsubsection{{Second order} Boltzmann Gibbs principle}
%
%In this section we follow the procedure given in the proof of Corollary 3.1 in \cite{ayala2018quantitative}. The general idea is to extract a factor $n^{-\gamma}$, for some $\gamma >0$, from the transition density of a pair of random walkers coupled to our dynamics.

\bl\label{secondBGprinciple}
 Let  $\caE^{(n,k)}(\phi, \eta)$ be given by \eqref{E}, then, for  every test function $\phi\in \hat S(\R^{d})$ there exists $C>0$ such that, for all $t\ge 0$ and $n\in \N$,
\be
 \E_n \left[ \( \int_0^t \caE^{(n,k)}(\phi, \eta(n^2s)) ds  \)^2 \right]  \le C\cdot   \frac{t^2}n.
\ee
\el

\bpr
%From \eqref{ciao1} we have that
%the error function
%is of the form:
%\beq\label{caEn}
%\caE_n(\phi^{(2)},\eta) &=& \frac{1}{n} \sum_{ \bix \in \Z^2} f_n(\bix)\cdot \D(\bix, \eta)  
%%\\ 
%%&=& \col{\frac{\sigma}{n} \sum_{x_1 \in \Z}  \left[ \sum_{ x_2 \in \Z} p(x_2-x_1) \(x_2-x_1\)^2  \D({x_2}, \eta)  + O(n^{-1}) \right] \(\phi\myprime\(\tfrac{x_1}n\)\)^2  \D({x_1}, \eta) } \nn 
%\eeq
%with
%\beq\label{ff}
%f_n(\bix)= E_n(\phi^{(2)})({\mathbf x}) +\tfrac 1 n\cdot   \phi\(\tfrac{x_1}n\) \psi^{(n)}\(\tfrac{x_2}n\)
%\eeq
%where
%\beq\label{splitgentwoparticlesSecondterm}
%E_n(\phi^{(2)})({\mathbf x})
%&=& \sigma n^2 \left[\phi_n({x_2})  p(x_2-x_1)  \(\phi\(\tfrac{x_2}n\)- \phi\(\tfrac {x_1}n\)\) \right] \nn \\
%&+& \sigma n^2 \left[\phi_n({x_1})  p(x_1-x_2)\(\phi\(\tfrac{x_2}n\)- \phi\(\tfrac {x_1}n\)\)\right] \nn \\
%&=& \sigma n^2  p(x_2-x_1) \(\phi(\tfrac{x_2}n) -\phi(\tfrac{x_1}n) \)^2 \nn \\
%&=& \sigma p(x_2-x_1) \(x_2-x_1\)^2 \(\phi\myprime(\tfrac{x_1}n)+\phi\myprime\myprime(\tfrac{x_1}n) \cdot \tfrac{x_2-x_1} n+\ldots\)^2 
%%&=& \sigma p(x_2-x_1) \(x_2-x_1\)^2 \(\sum_{m=1}^\infty \tfrac 1 {m!} \cdot \tfrac{d^m\phi}{dx^m}(\tfrac{x_1}n)\cdot \(\tfrac{x_2-x_1}n\)^{m-1}\)^2 
%\eeq
%where in the second equality we used the symmetry of $p(\cdot)$, and in the third we Taylor expanded around $x_1$ inside the braces. 
Using the fact  that $\phi$ is bounded and that $p(\cdot)$ has finite range we can conclude that there exists an $M>0$ such that
\beq\label{M}
\sup_n \sup_{\bix \in \Z^{kd}} |\Psi_n(\bix)|\le M.
\eeq
%Then, from \eqref{ciao1} we have
%\beq\label{caEn}
%\caE_n(\phi^{(2)},\eta) &=& \frac{1}{n} \sum_{ x_1, x_2 \in \Z} \left[ \sigma p(x_2-x_1) \(x_2-x_1\)^2 \(\phi\myprime\(\tfrac{x_1}n\)\)^2{+ O(n^{-1})}\right]  \D((x_1,x_2), \eta)  \nn 
%%\\ 
%%&=& \col{\frac{\sigma}{n} \sum_{x_1 \in \Z}  \left[ \sum_{ x_2 \in \Z} p(x_2-x_1) \(x_2-x_1\)^2  \D({x_2}, \eta)  + O(n^{-1}) \right] \(\phi\myprime\(\tfrac{x_1}n\)\)^2  \D({x_1}, \eta) } \nn 
%\eeq
%%\col{where the last identity follows from the fact that $\D((x_1,x_2),\eta)=\D(x_1,\eta)\cdot \D(x_2,\eta)$ for $x_1\neq x_2$ and from the assumption $p(0)=0$.}
%%\\
%%\col{is this used?}
%%\\
\noindent
We recall here that the duality function is parametrized by the density parameter $\rho$, i.e. $D(\cdot,\cdot)=D_\rho(\cdot,\cdot)$ and that  $\{D_\rho(\xi,\cdot), \xi\in \Omega\}$ is a family of products of polynomials that are orthogonal with respect to the reversible measure $\nu_\rho$. From the stationarity of $\nu_\rho$ we have
\beq\label{statnBG2part}
&&\E_n \left[ \( \int_0^t \caE^{(n,k)}(\phi, \eta(n^2s)) ds  \)^2 \right] =  \int_0^t \int_0^t \E_{n} \left[ \caE^{(n,k)}(\phi, \eta_{n^2s}) \caE^{(n,k)}(\phi, \eta_{n^2u} ) \right]du \, ds   \nn \\
&&\hskip3.5cm=2 \int_0^t \int_0^s \int   \E_{\eta} \left[ \caE^{(n,k)}(\phi, \eta_{n^2(s-u)}) \right] \caE^{(n,k)}(\phi, \eta )   \nu_\rho ( d \eta) du \, ds.
\eeq
The fact that we can exchange expectations and integral is a consequence of Proposition \ref{mombndXk} in  Section \ref{repsec}, which does not use any results of the current section. \\

\noindent
Let us denote by $V_n(\phi)$ the integrand in \eqref{statnBG2part},  then, using \eqref{ort},  we have
\beq\label{vfunctwopart}
V_n(\phi)&=&   \frac{1}{n^{kd}} \sum_{\bix, \biy\in \Z^{kd}} \Psi_n(\bix) \Psi_n(\biy) \cdot \pii(\bix) \pii(\biy)\cdot  \int   \E_{\eta} \left[D_\rho(\xi(\bix), \eta_{n^2(s-u)})  \right]   D_\rho(\xi(\biy), \eta)    \nu_\rho ( d \eta)\nn \\
&=&   \frac{1}{n^{kd}} \sum_{\bix, \biy\in  \Z^{kd}} \Psi_n(\bix) \Psi_n(\biy) \cdot \pii(\bix) \pii(\biy) 
\cdot  \frac 1{\mu_\rho(\xi(\biy))} \cdot p_{n^2(s-u)}(\xi(\bix), \xi(\biy))\nn\\
&=&  \frac{c}{n^{kd}}\sum_{\bix\in \Z^{kd}} \Psi_n(\bix)  \cdot \pii(\bix) \sum_{\biy\in \Z^{kd}}\frac 1 { N(\xi(\biy))}\cdot \Psi_n(\biy)
 \cdot p_{n^2(s-u)}(\xi(\bix), \xi(\biy))\nn\\
 &\le& \frac{cM}{n^{kd}} \sum_{\bix\in \Z^{kd}} \mid \Psi_n(\bix) \mid \cdot \pii(\bix) \sum_{\biy\in \Z^{kd}} \frac 1 {N(\xi(\biy))}
 \cdot p_{n^2(s-u)}(\xi(\bix), \xi(\biy))\nn\\
  &=& \frac{cM}{n^{kd}} \sum_{\bix\in  \Z^{kd}} \mid \Psi_n(\bix) \mid  \cdot \pii(\bix) \sum_{\xi'\in\Omega_k} 
 \cdot p_{n^2(s-u)}(\xi(\bix), \xi')\nn\\
   &\le& \frac{c'M}{n^{kd}} \sum_{\bix\in  \Z^{kd}} \mid \Psi_n(\bix) \mid 
 \eeq
where  we used \eqref{ort} in the second identity, \eqref{pf} and \eqref{La} in the third identity (with $c=c(k,\rho)$) and \eqref{M} in the fourth line. From \eqref{PSI} we have
\beq
\frac{1}{n^{kd}} \sum_{\bix\in \Z^{kd}} \mid \Psi_n(\bix) \mid  &\le& {\frac \alpha{n^{kd+1}}\sum_{\bix \in \Z^{kd}}  \mid \psi_{n} \mid \(\tfrac{x_k}n\)\cdot \prod_{\ell=1}^{k-1} \mid \phi(\tfrac{x_\ell}n)} \mid \label{ciao3}\\
&+& \frac{\sigma(k-1)}{2n^{kd}} \sum_{\bix\in \Z^{kd}}   \prod_{\ell=3}^{k} \mid \phi(\tfrac{x_\ell}n) \mid \cdot
 p(x_2-x_{1})n^2\(\phi(\tfrac{x_{2}}n)-\phi(\tfrac{x_1}n)\)^2.\nn
 \eeq
Using \eqref{QUO} we have that the first term in the r.h.s. of \eqref{ciao3} is bounded by a constant times $n^{-1}$.
For what concerns   the second term, we have:
 \beq
 && \frac{\sigma(k-1)}{2n^{(k-2)d}}\(\prod_{\ell=3}^k\sum_{x_\ell\in \Zd} \phi( \tfrac{x_\ell}n) \)\cdot \frac 1 {n^{2d}}\sum_{x_1,x_2\in \Z^{d}}
 p(x_2-x_1) n^2\(\phi\(\tfrac{x_{2}}n\)-\phi\(\tfrac{x_1}n\)\)^2\nn\\
 &&\le   \frac c {n^{2d}}\sum_{x_1,x_2\in \Z^{d}}
 p(x_2-x_1) n^2\(\phi\(\tfrac{x_{2}}n\)-\phi\(\tfrac{x_1}n\)\)^2\nn
 \eeq
 Now, from  the Taylor expansion \eqref{taylor} we know  that there exists a sequence of functions
 where, using   the fact that the range of $p(\cdot)$ is $\caR=[-R,R]^d$, and the Taylor expansion \eqref{taylor} we have that there exists a smooth function $\widetilde \psi\in S(\R^d)$ such that, for all $x\in \Zd$,
 \beq
 \sup_{r\in \R} \{n^2\(\phi\(\tfrac{x+r}n\)-\phi\(\tfrac{x}n\)\)^2 \}\le \widetilde \psi\(\tfrac x n\)
 \eeq
 as a consequence we obtain the upper bound
 \beq
 &&\frac 1 {n^{2d}}\sum_{x_1,x_2\in \Z^{d}}
 p(x_2-x_1)n^2\(\phi(\tfrac{x_{2}}n)-\phi(\tfrac{x_1}n)\)^2
%&=& \frac 1 {n^{2d}}\sum_{\bix\in \Z^2}  p(x_2-x_1) \(x_2-x_1\)^2 \(\phi\myprime(\tfrac{x_1}n)+\phi\myprime\myprime(\tfrac{x_1}n) \cdot \tfrac{x_2-x_1} n+\ldots\)^2 \label{st}\nn\\
=\frac 1 {n^{2d}}\sum_{x\in \Z^{d}}\sum_{r\in \caR}
 p(r)n^2\(\phi\(\tfrac{x+r}n\)-\phi\(\tfrac{x}n\)\)^2\nn\\
&&\le \frac{1}{n^{2d}} \sum_{r\in \caR} \sum_{x\in \Zd}  p(r)  \cdot  \widetilde \psi\(\tfrac xn\)\le \frac{c}{n^{2d}}  \sum_{x\in \Zd}   \widetilde \psi\(\tfrac{x}n\)\le \frac{c'}{n^d}
\eeq
where the inequality holds for a suitable $c'>0$. %Moreover, given $x_1 \in \Z$ there exists $R > 0$ such that
%\be
%\sup_{a, b \in \Z} p(a-b) \(a-b\)^2 \leq R
%\ee
%hence
In conclusion we have that there exists a constant $C>0$ such that
\beq\label{vfunctwopartesti2}
V_n(\phi) &\leq&   \frac{C}{n} 
\eeq
from which the statement follows.
\epr
\vskip.3cm
\noindent
As a consequence  of Lemma \ref{secondBGprinciple} we can  close the drift term, i.e.
\beq\label{closeddrifttwo}
\loc \caX^{(n,k)}(\Phi,\eta) &=&\alpha k \cdot\tfrac {\chi}2 \cdot \caY^{(n,k)}(\phi^{(k-1)}\otimes \Delta\phi,\eta)+ \caE^{(n,k)}(\phi,\eta) \nn \\
&=& \alpha k \cdot\tfrac {\chi}2 \cdot \caY^{(n,k)}(\phi^{(k-1)}\otimes \Delta\phi,\eta) + O(n^{-1}).
\eeq

\subsection{Closing the equation for the carr\'e-du-champ}\label{CdC}
In this section we will show that  the integrand in the RHS of equation \eqref{Ntdef} can be expressed in terms of the $(k-1)$th-order fluctuation field $\caX^{(n,k-1)}$. To achieve this we consider the expression  for the carr\'e-du-champ given by \eqref{Eqsquareduchamp} in the Appendix. For the case of our $k$th-order fluctuation field this becomes
\be\label{Ga}
n^2 \Gamma \caX^{(n,k)}(\Phi,\eta) = \frac{1}{n^d} \sum_{x \in \Zd } \sum_{r \in \caR} p(r) \eta_x ( \alpha + \sigma \eta_{x+r})  \left[n^{d/2+1} \(\caX^{(n,k)}(\Phi,\eta^{x,x+r})-\caX^{(n,k)}(\Phi,\eta)\)  \right]^2.
\ee
Notice that here we multiplied  by a factor $n^{d/2+1}$  the squared term in order  to cancel the $n^2$ in front of the carr\'e-du-champ and get a general factor $n^{-d}$ in front of the sum.\\

\noindent
In the next section we find some recursion relations for duality polynomials. The main application of these relations consists in allowing us to rewrite  any polynomial depending on $\eta^{x,x+r}$ in terms of polynomials depending on the unmodified $\eta$. 

\subsubsection{Recursion relation for duality polynomials}

In this section  we obtain   a recurrence relation for the single-site orthogonal polynomials. Before giving the result it is convenient to 
 summarize  the expression for the self-duality  generating function by defining the function
\beq\label{genfun}
f_\sigma(t,n):=\sum_{m=0}^{\infty}  \di(m,n)\cdot{t^m}
\eeq
then $f_\sigma$ can be written in the form
\be\label{hsip}
f_\sigma(t,n) =e_\sigma(t)\cdot h_\sigma(t)^n, \qquad  h_\sigma(t) =  \frac{ 1 - c_{-\sigma} {t}}{1-\sigma t}, \qquad  e_\sigma(t)=
\left\{
\begin{array}{lll}
e^{-t}& \text{if} & \sigma=0\\
(1+\sigma t)^{\sigma \alpha} & \text{if} & \sigma=\pm 1
\end{array}
\right.
\ee
with $\cs$ given by \eqref{cs}.
%\be\label{cs}
 %\cs=
%\left\{
%\begin{array}{lll}
%1/ \rho & \text{if} &\sigma=0\\
%{\alpha}/{\rho} & \text{if}& \sigma=1\\
%(\alpha+\rho)/{\rho} & \text{if} & \sigma=-1
%\end{array}
%\right..
%\ee 
%\col{$\cs$ should be defined before the main theorem}
Then we define the functions $g_\sigma, \tilde g_\sigma: \N \to \R$   given by
\beq
&&g_\sigma(m) :=\frac 1 {m!} 
\left.\frac{d^{m}}{dt^{m}}h_\sigma(t) \right|_{t=0} 
\qquad \text{and} \qquad
\tilde{g}_\sigma(m) :=\frac 1 {m!} 
\left. \frac{d^{m}}{dt^{m}}\frac{1}{h_\sigma(t)} \right|_{t=0}  \qquad \text{for } m\ge 1\nn\\
&&\nn\\
&&\text{and} \qquad g_\sigma(0)=\tilde g_\sigma(0):=1\label{g0}
\eeq
that are exactly computable:
\be\nonumber
g_\sigma(m)=
\left\{
\begin{array}{ll}
-\frac 1 \rho\cdot  \, \mathbf 1_{m=1} & \sigma=0\\
&\\
-\frac \alpha \rho & \sigma=+1\\
&\\
-\frac{\alpha+\rho}\rho\cdot \rho (-1)^{m-1} & \sigma=-1
\end{array}
\right.
\qquad \tilde g_\sigma(m)=
\left\{
\begin{array}{ll}
\(\frac 1 \rho\)^m  & \sigma=0\\
&\\
\(\frac \alpha \rho\) \cdot \(\frac{\alpha+\rho}{\rho}\)^{m-1}  & \sigma=+1\\
&\\
\frac{\alpha+\rho}{\rho} \cdot \(\frac \alpha \rho\)^{m-1}& \sigma=-1
\end{array}
\right.
\ee
for $m\ge 1$,
that can be rewritten as
\beq\label{g}
g_\sigma(m)= -\cs \cdot \sigma^{m-1}\qquad \text{and}\qquad \tilde g_\sigma(m)=\cs\cdot c_{-\sigma}^{m-1}  \qquad \text{for }m\ge 1
\eeq
notice, in particular, that
\be\label{g1}
\tilde g_\sigma(1)=\cs=-g_\sigma(1).
\ee
 We have the following result.
\bt\label{propcontrolbeta}
For any $m,n \in \N$ we have
\be\label{propcontrolbetaeq}
\di(m,n+1) = \sum_{j=0}^{m} g(m-j) \cdot \di(j,n)    \\
\ee
and
\be\label{propcontrolbetaminuseq}
\di(m,n-1) = \sum_{j=0}^{m}  \tilde{g}(m-j)\cdot  \di(j,n)     \\
\ee
with $g,\tilde g: \N\to \R$  as in \eqref{g0}-\eqref{g}.
\et

\bpr
From \eqref{hsip} we have that 
\be\label{decompfh}
f(t,n+1) = f(t,n) h(t) 
\ee
then, from the generating function definition \eqref{genfun}, we  deduce that
\be\label{diffMeixfrel}
\di(m,n) =  \frac 1 {m!}\cdot \left. \frac{d^{m}}{dt^{m}} f(t,n) \right|_{t=0}
\ee
hence, the recurrence relation \eqref{decompfh} and an application of Leibniz product rule for differentiation in the RHS above give
\beq\label{Mlbplusone}
\di(m,n+1) &=&   \frac 1 {m!}\cdot \sum_{j=0}^{m} \binom{m}{j} \left. \frac{d^{j}}{dt^{j}} f(t,n)\right|_{t=0} \left. \cdot \frac{d^{m-j}}{dt^{m-j}}h(t) \right|_{t=0}\nn \\
&=&   \frac 1 {m!}\cdot\sum_{j=0}^{m} \binom{m}{j} j!\cdot  \di(j,n) \left. \cdot  \frac{d^{m-j}}{dt^{m-j}}h(t) \right|_{t=0} \nn\\
&=&\sum_{j=0}^{m} \frac 1 {(m-j)!}    \cdot \frac{d^{m-j}}{dt^{m-j}}h(t) \bigg|_{t=0} \cdot \di (j,n)\nn\\
&=&\sum_{j=0}^{m}  g(m-j) \cdot \di(j,n)\nn
\eeq
where in the second equality we used \eqref{diffMeixfrel}. This concludes the proof of \eqref{propcontrolbetaeq}.  Equation \eqref{propcontrolbetaminuseq}
%%
%%\noindent
%%
%%%\frac{\Gamma(\alpha + j)}{\Gamma(\alpha+c)} \left. \frac{d^{c-s}}{dt^{c-s}}h(t) \right|_{t=0}  \nn
%%Moreover, for any natural number $m\neq 0$ a simple iterated differentiation over $h$ gives the expression
%%\be
%%\left. \frac{d^{m}}{dt^{m}}h(t) \right|_{t=0} = - m! \( \frac{\alpha}{\rho} \) \nn
%%\ee
%%hence $g: \N \to \R$ can be re-expressed as 
%%\be\label{gdef}
%%g(m) := \begin{dcases*}
%%-m! \( \frac{\alpha}{\rho} \) & if $m>0$\\
%%1 & if $m=0$
%%\end{dcases*}
%%\ee
%%hence by \eqref{Mlbplusone} we have
%%\beq\label{Mlbplusonegexp}
%%d(c,a+1) \pi(c) &=& M(c,a+1) \nn \\
%%&=&  \sum_{s=0}^{c} {c \choose s} M(s,a) g(c-s) \nn \\
%%&=&  \sum_{s=0}^{c} {c \choose s} d(s,a)\pi(s) g(c-s) \nn \\
%%\eeq
%%and the proof is completed.
%\epr
%
%\bpr[Proposition \ref{propcontrolbetaminus}]\\
can be proved from the same reasoning,  with the difference that we now have the inverse relation
\be\label{decompfhinv}
f(t,n-1) = f(t,n) \cdot \frac{1}{h(t)}. 
\ee
This change results, after the application of Leibniz  rule, in the relation
\beq\label{Mlbminusone}
\di(m,n-1) &=&  \frac 1 {m!}\cdot \sum_{j=0}^{m} {m \choose j} j! \cdot \di(j,n) \left. \cdot \frac{d^{m-j}}{dt^{m-j}}\frac{1}{h(t)}\right|_{t=0} \nn \\
&=&  \sum_{j=0}^{m}  \tilde g(m-j) \cdot \di(j,n)\nn
\eeq
that concludes the proof.
\epr

\subsubsection{Controlling the moments of the fields}\label{repsec}
The objective of this section is to take advantage of the ergodic properties of our process to introduce a result that will allow us to make multiple replacements, in the appropriate sense, inside the expression of the carr\'e-du-champ given in \eqref{Ga}. Let us start first with a uniform estimate for moments of the fields $\caX^{(n,l)}(\Phi,\eta)$.
\bp\label{mombndXk}
Let $l, m \in \N$ then we have
\be
\sup_{n \in \N} \, \E_{\nu_\rho} \left[ \caX^{(n,l)}(\Phi,\eta )^m \right] \leq C(\rho,\phi)
\ee
\ep
\bpr
As claimed in the statement of the proposition, this result holds for any finite natural number $m$. Nevertheless for simplicity we will only show how to obtain the estimates for $m \in \{2,4 \}$ (which indeed are the only two uses that we make of this result). Let us start with the simplest non-trivial case,  $m=2$,  for which the result comes directly from orthogonality
\beq\label{unifbdkl}
\E_{\nu_\rho} \left[ \caX^{(n,l)}(\Phi,\eta)^2 \right] &=& n^{-ld} \sum_{\xi , \xi' \in\Omega_l}  \Phi_n(\xi) \Phi_n(\xi') \La(\xi) \La(\xi')  \E_{\nu_\rho} \left[ D(\xi,\eta) D(\xi',\eta) \right]  \\
\label{unifbdkl2}
&=& n^{-ld} \sum_{\xi \in\Omega_l}  \Phi_n(\xi)^2  \La(\xi)^2  \frac{1}{\mu_\rho(\xi)}  \\
&\leq& K \cdot n^{-ld} \sum_{\xi \in\Omega_l}  \Phi_n(\xi)^2  < \infty
\eeq
where in the second line we used \eqref{Ort} and $K$ is given by 
\be
K = \sup_{\xi \in \Omega_k}  \, \frac{\La(\xi)^2}{\mu_\rho(\xi)}. \nn
\ee
Notice that the previous estimate was possible due the fact that orthogonality, in the form of expression \eqref{Ort}, allowed us to reduce the summation in the RHS of \eqref{unifbdkl} from a $2ld$ dimensional sum to an $ld$ dimensional sum in \eqref{unifbdkl2}.\\
\vskip0.2cm
\noindent
For the case $m=4$ we have
\beq\label{unifbdkl4}
\E_{\nu_\rho} \left[ \caX^{(n,l)}(\Phi,\eta)^4 \right] &=& n^{-2ld} \sum_{\xi^{(j)}\in \Omega_l}  \prod_{j=1}^4 \Phi_n(\xi^{(j)}) \cdot \La(\xi^{(j)}) \cdot \nn \\
&& \E_{\nu_\rho} \left[ D(\xi^{(1)},\eta) D(\xi^{(2)},\eta)D(\xi^{(3)},\eta) D(\xi^{(4)},\eta) \right] 
\eeq
For this case the sum in the RHS of \eqref{unifbdkl4} is $4ld$-dimensional. Given the factor $n^{-2ld}$ in front of the RHS, in order to obtain a uniform estimate, we would like this summation to be $2ld$ dimensional instead. In order to see that this is indeed the case, we analyze the non-zero contribution coming from
\be
\E_{\nu_\rho} \left[ D(\xi^{(1)},\eta) D(\xi^{(2)},\eta)D(\xi^{(3)},\eta) D(\xi^{(4)},\eta) \right]. \nn
\ee
By the product nature of the measure $\nu_\rho$ and the duality polynomials we have
\be\label{prod4dual}
\E_{\nu_\rho} \left[ D(\xi^{(1)},\eta) D(\xi^{(2)},\eta)D(\xi^{(3)},\eta) D(\xi^{(4)},\eta) \right] = \prod_{x \in \Zd} \E_{\nu_\rho} \left[ d(\xi_x^{(1)},\eta) d(\xi_x^{(2)},\eta)d(\xi_x^{(3)},\eta) d(\xi_x^{(4)},\eta) \right].
\ee
Notice that for every $x$ for which $\xi_x^{(j)}=0$  for all $j \in \{1,2,3,4\}$, the corresponding contribution in the RHS of \eqref{prod4dual} is equal to $1$ and therefore negligible. This is precisely the reason why the summation in the RHS of \eqref{unifbdkl4} is at most $4ld$-dimensional. We have indeed that  the maximum number of $x \in \Zd$  contributing to the product in the RHS of \eqref{prod4dual} is at most $4l$, i.e. one for each of the $4l$ particles that all the $\xi^{(j)}$ have in total. In reality we can see that there are less $x$s giving a non-zero contribution. In order to see is, consider  an $x \in \Zd$ such that there exists a unique $j \in \{1,2,3,4\}$ for which $\xi_x^{(j)}\neq 0$. In this case,  because of the zero mean of the single-site duality function we have
\be
\E_{\nu_\rho} \left[ d(\xi_x^{(1)},\eta) d(\xi_x^{(2)},\eta)d(\xi_x^{(3)},\eta) d(\xi_x^{(4)},\eta) \right] = 0
\ee
this means that whenever  $x \in \Zd$ is such that there exists a $j \in \{1,2,3,4\}$ for which $\xi_x^{(j)}\neq 0$ there must be another $j' \in \{1,2,3,4\}$ for which $\xi_x^{(j')}\neq 0$. In other words we only have a possibility of $2l$ particles to distribute freely, and hence the summation in the RHS of \eqref{unifbdkl4} is at most $2ld$-dimensional.
\epr

\bp\label{metarep}
Let $f: \Rd \to \R$ be a test function, and $\{ M_n: \Omega \times \R  \to \R : n \in \N \}$ be  a sequence of uniformly bounded cylindrical functions of the form
\be\label{requirMn}
M_n(\eta,x) = f(x/n) \prod_{j \in \N} d(b_j,\eta_x) 
\ee
where only a finite number of $b_j$ are different from zero.  Let also $\{a_n : n \in \N\}$ be a sequence of real numbers converging to $0$, we then have 
\beq
\lim_{n\to \infty} \E_n\left[\ \(\int_0^t \frac{a_n}{n^d} \sum_{x \in \Zd } \sum_{r \in \caR} p(r) \eta_x(n^2 s) ( \alpha + \sigma \eta_{x+r}(n^2 s))   M_n(\eta_x(n^2s)) \cdot
  \caX^{(n,l)}(\Phi,\eta(n^2s))^m\; ds \)^2 \right]=0 \nn 
\eeq                       
for all $l \in \{1,2,\dots,k-1\}$, and  $m \in \N$.
\ep

\bpr
By Cauchy-Schwartz we have
\beq
&& \E_n\left[\ \(\int_0^t \frac{a_n}{n^d} \sum_{x \in \Zd } \sum_{r \in \caR} p(r) \eta_x(n^2 s) ( \alpha + \sigma \eta_{x+r}(n^2 s))  \cdot M_n(\eta_x(n^2s)) \cdot
 \caX^{(n,l)}(\Phi,\eta(n^2s))^m\; ds \)^2 \right]  \nn \\
&\leq&  \frac{a_n^2 t}{n^{2d}}  \int_0^t \E_{n} \left[  \caX^{(n,l)}(\Phi,\eta(n^2s))^{2m} \cdot \( \sum_{x \in \Zd } \sum_{r \in \caR} p(r) \eta_x(n^2 s) ( \alpha + \sigma \eta_{x+r}(n^2 s)) \cdot M_n(\eta_x(n^2s)) \)^2 \right] \; ds \nn \\
&=&  \frac{a_n^2 t^2}{n^{2d}} \E_{n} \left[ \caX^{(n,l)}(\Phi,\eta)^{2m} \cdot \( \sum_{x \in \Zd } \sum_{r \in \caR} p(r) \eta_x ( \alpha + \sigma \eta_{x+r})  \cdot M_n(\eta_x) \cdot
   \)^2 \right]  \nn \\
&=&  \frac{a_n^2 t^2}{n^{2d}} \sum_{x,y \in \Zd } \sum_{r_1, r_2 \in \caR} p(r_1) \cdot p(r_2) \cdot \E_{n} \left[ M_n(\eta_x)  \cdot M_n(\eta_y) \cdot \caX^{(n,l)}(\Phi,\eta)^{2m}  \right]  \nn \\
&\le& \frac{a_n^2 t^2}{n^{2d}} \sum_{x,y \in \Zd } \sum_{r_1, r_2 \in \caR} p(r_1) \cdot p(r_2) \cdot \sqrt{ \E_{n} \left[ M_n(\eta_x)^2  \cdot M_n(\eta_y)^2 \right] } \cdot \sqrt{\E_{n} \left[ \cdot \caX^{(n,l)}(\Phi,\eta)^{4m} \right]}  \nn \\
&\le& K t^2 a_n^2
\eeq
where in the last line we used Proposition \ref{mombndXk}, the boundedness of the single-site duality polynomials $d(b_j,\eta_x)$ and the smoothness of $f$ in the representation \eqref{requirMn}. The result then follows from the convergence  $a_n \to 0.$
\epr

\subsubsection{The gradient of the fluctuation fields}
Our goal for this section is to rewrite the square inside the RHS of \eqref{Ga} in terms of lower order fluctuation fields. We will see that this can  be expressed, in agreement with \eqref{QuadvarmartingXk} , only in terms of the field of order $k-1$. Let us  then denote by $\nabla_{d}^{i,i+r}$ the $d$-dimensional gradient
\be\label{defgradkthfield}
\nabla_{d}^{i,i+r} \caX^{(n,k)}(\Phi,\eta) = n^{d/2+1}\( \caX^{(n,k)}(\Phi,\eta^{i,i+r})-\caX^{(n,k)}(\Phi,\eta) \).
\ee 
\noindent
Notice that, by linearity of the $k$-th order field, we have
\beq\label{compgradkthfield1}
\nabla^{i,j}\caX^{(n,k)} (\Phi,\eta)&:=& n^{-\tfrac{(k-1)d}{2}+1} \sum_{ \xi \in \Omega_k} \Phi_n(\xi) \left[ \Di(\xi, \eta^{i,j}) -\Di(\xi, \eta) \right]
\eeq
with $\Di(\cdot,\cdot)$ as in \eqref{prod}.
We define now, for $i,j\in\Zd$, $\ell\le k$, the auxiliary field
\beq\label{auxfieldZ}
{\cal Z}^{(n,k,\ell)}_{i,j}(\Phi,\eta) :=n^{-kd/2} \sum_{\xi\in\Omega_k} \mathbf 1_{\xi_i+\xi_j=\ell}\cdot \Phi_n(\xi)\Di(\xi,\eta)
\eeq
then we have the following formula for the gradient of the fluctuation field.
\bp\label{PropformGrad}
\beq\label{sum}
\nabla^{i,j}\caX^{(n,k)} (\Phi,\eta)=\sum_{s=1}^k n^{-\tfrac{(s-1)d}2}\cdot \sum_{m=1}^{s} n\(\phi(\tfrac j n)^{m} -\phi(\tfrac in)^{m}\) \cdot \gi(m)\cdot {\cal Z}^{(n,k-m,s-m)}_{i,j}(\phi,\eta-\delta_i)  \nn
\eeq
\ep
\bpr
Using the product nature of the  polynomials $\Di(\cdot,\eta)$  and of $\Phi_n(\cdot)$ we get
\beq
\nabla^{i,j}\caX^{(n,k)} (\Phi,\eta)=n\sum_{s=1}^k n^{-\tfrac{(s-1)d}2}\cdot {\cal Z}^{(n,k-s,0)}_{i,j}(\phi,\eta)  \cdot \sum_{a=0}^s Y^{(n,a,s-a)}_{i,j}(\phi,\eta) 
\eeq
and
\beq
&&Y^{(n,a,b)}_{i,j}(\phi,\eta):=\phi(\tfrac in)^a\phi(\tfrac jn)^b
\left\{\di(a,\eta_i-1)\di(b,\eta_j+1)-\di(a,\eta_i)\di(b,\eta_j)\right\}\nn\\
&&\hskip1cm=\phi(\tfrac in)^a\phi(\tfrac jn)^b
\left\{\di(a,\eta_i-1)\left[\di(b,\eta_j+1)-\di(b,\eta_j)\right]+\di(b,\eta_j)\left[\di(a,\eta_i-1)-\di(a,\eta_i)\right]\right\}\nn
\eeq
hence, using \eqref{propcontrolbetaeq}
 we get
\beq
  \sum_{a=0}^s Y^{(n,a,s-a)}_{i,j}(\phi,\eta)&=&
\sum_{a=0}^{s-1} \phi(\tfrac in)^a\phi(\tfrac jn)^{s-a}
\di(a,\eta_i-1)\left[\di(s-a,\eta_j+1)-\di(s-a,\eta_j)\right]\nn\\
&-&\sum_{b=0}^{s-1} \phi(\tfrac in)^{s-b}\phi(\tfrac jn)^{b}
\di(b,\eta_j)\left[\di(s-b,\eta_i)-\di(s-b,\eta_i-1)\right]\nn\\
&=&\sum_{a=0}^{s-1} \sum_{\kappa=0}^{s-a-1}\phi(\tfrac in)^a\phi(\tfrac jn)^{s-a} \cdot \gi(s-a-\kappa)\cdot 
\di(a,\eta_i-1) \di(\kappa,\eta_j)\nn\\
&-&\sum_{b=0}^{s-1} \sum_{m=0}^{s-b-1}\phi(\tfrac in)^{s-b}\phi(\tfrac jn)^{b} \cdot  \gi(s-b-m)\cdot
 \di(m,\eta_i-1)\di(b,\eta_j)\nn
\eeq
now, calling $b=\kappa$ and $m=a$ we get
\beq
&&\sum_{a=0}^{s-1} \sum_{\kappa=0}^{s-a-1}  \phi(\tfrac in)^a\phi(\tfrac jn)^{s-a} \cdot \gi(s-a-\kappa)\cdot
\di(a,\eta_i-1) \di(\kappa,\eta_j)\nn\\
&&-\sum_{\kappa=0}^{s-1} \sum_{a=0}^{s-\kappa-1} \phi(\tfrac in)^{s-\kappa}\phi(\tfrac jn)^{\kappa} \cdot\gi(s-a-\kappa)\cdot
 \di(a,\eta_i-1)\di(\kappa,\eta_j)\nn\\
 &&=\sum_{a=0}^{s-1} \sum_{\ell=a}^{s-1}  \(\phi(\tfrac in)^a\phi(\tfrac jn)^{s-a} -\phi(\tfrac in)^{s+a-\ell}\phi(\tfrac jn)^{\ell-a} \)\cdot \gi(s-\ell)\cdot
\di(a,\eta_i-1) \di(\ell-a,\eta_j)\nn\\
 &&=\sum_{a=0}^{s-1} \sum_{\ell=a}^{s-1}  \(\phi(\tfrac jn)^{s-\ell} -\phi(\tfrac in)^{s-\ell}\)\phi(\tfrac in)^a\phi(\tfrac jn)^{\ell-a}\cdot \gi(s-\ell)\cdot
\di(a,\eta_i-1) \di(\ell-a,\eta_j)\nn\\
 &&=\sum_{\ell=0}^{s-1} \(\phi(\tfrac jn)^{s-\ell} -\phi(\tfrac in)^{s-\ell}\) \cdot \gi(s-\ell)\sum_{a=0}^{\ell}  \phi(\tfrac in)^a\phi(\tfrac jn)^{\ell-a}\cdot
\di(a,\eta_i-1) \di(\ell-a,\eta_j)\nn
\eeq
where the first  identity follows from the change of variable $\ell=\kappa+a$.
Then
\beq
&&{\cal Z}^{(n,k-s)}_{i,j}(\Phi,\eta)  \cdot \sum_{a=0}^s Y^{(n,a,s-a)}_{i,j}(\phi,\eta)\nn\\
&&= \sum_{\ell=0}^{s-1} \(\phi(\tfrac jn)^{s-\ell} -\phi(\tfrac in)^{s-\ell}\) \cdot \gi(s-\ell)\cdot {\cal Z}^{(n,k-s,0)}_{i,j}(\Phi,\eta)  \cdot \sum_{a=0}^{\ell}  \phi(\tfrac in)^a\phi(\tfrac jn)^{\ell-a}\cdot
\di(a,\eta_i-1) \di(\ell-a,\eta_j)\nn\\
&&= \sum_{\ell=0}^{s-1} \(\phi(\tfrac jn)^{s-\ell} -\phi(\tfrac in)^{s-\ell}\) \cdot \gi(s-\ell)\cdot {\cal Z}^{(n,k-(s-\ell),\ell)}_{i,j}(\Phi,\eta-\delta_i)  
\eeq
then
\beq
\nabla^{i,j}\caX^{(n,k)} (\Phi,\eta)&=&n\sum_{s=1}^k n^{-\tfrac{(s-1)d}2}\cdot \sum_{\ell=0}^{s-1} \(\phi(\tfrac jn)^{s-\ell} -\phi(\tfrac in)^{s-\ell}\) \cdot \gi(s-\ell)\cdot {\cal Z}^{(n,k-(s-\ell),\ell)}_{i,j}(\Phi,\eta-\delta_i)  \nn\\
&=&\sum_{s=1}^k n^{-\tfrac{(s-1)d}2}\cdot \sum_{m=1}^{s} n\(\phi(\tfrac jn)^{m} -\phi(\tfrac in)^{m}\) \cdot \gi(m)\cdot {\cal Z}^{(n,k-m,s-m)}_{i,j}(\Phi,\eta-\delta_i).  \nn
\eeq
This concludes the proof.
\epr
\vskip.3cm
\noindent
The advantage that Proposition \ref{PropformGrad} gives us is that we now have an expression in terms of the auxiliary field \eqref{auxfieldZ}:
\beq\label{expstrat}
\nabla^{i,j}\caX^{(n,k)} (\Phi,\eta)&=&\sum_{s=1}^k n^{-\tfrac{(s-1)d}2}\cdot \sum_{m=1}^{s} n\(\phi(\tfrac jn)^{m} -\phi(\tfrac in)^{m}\) \cdot \gi(m)\cdot {\cal Z}^{(n,k-m,s-m)}_{i,j}(\Phi,\eta-\delta_i)  \nn \\
&=& n \(\phi(\tfrac jn) -\phi(\tfrac in)\) \cdot \gi(1)\cdot {\cal Z}^{(n,k-1,0)}_{i,j}(\Phi,\eta-\delta_i) \label{expstrat} \\
&+& \sum_{s=2}^k n^{-\tfrac{(s-1)d}2}\cdot \sum_{m=1}^{s} n\(\phi(\tfrac jn)^{m} -\phi(\tfrac in)^{m}\) \cdot \gi(m)\cdot {\cal Z}^{(n,k-m,s-m)}_{i,j}(\Phi,\eta-\delta_i).  \nn 
\eeq

\vskip.2cm
\noindent
Recall that we claimed that we are able to close the carr\'e-du-champ in an expression depending only on the field of order $k-1$. In order to achieve this  it remains to:
\ben
\item replace the first term in the RHS of \eqref{expstrat} by some expressions depending on the field of order $k-1$;
\item show that  the second term in the RHS of \eqref{expstrat} vanishes as $n \to \infty$.
\een 
We will achieve this in several steps, the first one being the proof of the following proposition.
\bp\label{Prop1}
For all $k\in \N$ we have
\beq
&&\lim_{n\to \infty} \E_n\left[\left(\int_0^t \frac{1}{n^d} \sum_{x \in \Zd } \sum_{r \in \caR} p(r) \eta_x(n^2 s) ( \alpha + \sigma \eta_{x+r}(n^2 s)) \right. \right. \nn \\
&&\hskip4cm\cdot \left. \left. \({\cal Z}^{(n,k,0)}_{x,x+r}(\Phi,\eta(n^2s)) -\caX^{(n,k)}(\Phi,\eta(n^2s))\)^2 \; ds \right)^2 \right]=0.
\eeq
\ep
\bpr
Notice that for any fixed $x$ we have
\beq
 \caX^{(n,k)}(\Phi,\eta(n^2s)) = \sum_{l=0}^k {\cal Z}^{(n,k,l)}_{x,x+r}(\Phi,\eta(n^2s))
\eeq
which implies
\beq
\({\cal Z}^{(n,k,0)}_{x,x+r}(\Phi,\eta(n^2s)) -\caX^{(n,k)}(\Phi,\eta(n^2s))\)^2 &=& \( \sum_{l=1}^k {\cal Z}^{(n,k,l)}_{x,x+r}(\Phi,\eta(n^2s)) \)^2 \nn \\
&\leq& k \sum_{l=1}^k {\cal Z}^{(n,k,l)}_{x,x+r}(\Phi,\eta(n^2s))^2. \nn
\eeq
Moreover, we can also estimate each ${\cal Z}^{(n,k,l)}_{x,x+r}(\Phi,\eta(n^2s))$ in terms of the coordinates field $\caY^{(n,k-l)}$ given by \eqref{defkorder} as follows:
\be
{\cal Z}^{(n,k,l)}_{x,x+r}(\Phi,\eta)^2 \le n^{-ld/2} \, M_n (\eta,l) \cdot \caY^{(n,k-l)}(\phi^{(k-l)},\eta)^2
\ee
where $M_n$ is made of terms of the form \eqref{requirMn}, i.e.
\be
M_n(\eta,l) = \sum_{\xi_x=0}^l \Phi(\xi_x \delta_x + ( l-\xi_x) \delta_{x+r}) \cdot d(\xi_x,\eta_x) \cdot d(l-\xi_x,\eta_{x+r}).
\ee
Thanks to Proposition \ref{metarep} we conclude the proof.
\epr

\vskip.2cm
\noindent
For what concerns the second step, let us denote by $\widehat { G}_{i,j}^{(n,k)} (\Phi,\eta)$ the second term in the RHS of \eqref{expstrat}, i.e.
\beq
\widehat { G}_{i,j}^{(n,k)} (\Phi,\eta):=\sum_{s=2}^k n^{-\tfrac{(s-1)d}2}\cdot \sum_{m=1}^{s} n\(\phi(\tfrac jn)^{m} -\phi(\tfrac in)^{m}\) \cdot \gi(m)\cdot {\cal Z}^{(n,k-m,s-m)}_{i,j}(\Phi,\eta-\delta_i)  \nn
\eeq
we have the following result  supporting our claim:
\bp\label{Prop2} Under the inductive hypothesis \ref{InducHyp} we have 
\beq
\lim_{n\to \infty} \E_n\left[\int_0^t\ \(\frac{1}{n^d} \sum_{x \in \Zd } \sum_{r \in \caR} p(r) \eta_x ( \alpha + \sigma \eta_{x+r})  \cdot   
  \widehat G_{x,x+r}^{(n,k)}(\phi,\eta(n^2 s))^2 \; ds \)^2\right]=0.
\eeq
\ep
\bpr
 After  expanding $\widehat G_{x,x+r}^{(n,k)}(\phi,\eta(n^2 s))^2$, the statement   follows from applying multiple times Propositions \ref{Prop1} and  \ref{metarep}.
\epr

\bp\label{Prop3}  Let
\be\label{error}
  G_{i,j}^{(n,k)}(\phi,\eta):= \nabla^{i,j} \caX^{(n,k)}(\Phi,\eta)  + \cs\,\langle j-i, \nabla\phi(\tfrac{i}n)\rangle \cdot \caX^{(n,k-1)}(\Phi,\eta) \nn
\ee
then, under the inductive hypothesis \ref{InducHyp}, we have
\beq
\lim_{n\to \infty} \E_n\left[\int_0^t\ \(\frac{1}{n^d} \sum_{x \in \Zd } \sum_{r \in \caR} p(r) \eta_x ( \alpha + \sigma \eta_{x+r})  \cdot  G_{x,x+r}^{(n,k)}(\phi,\eta(n^2 s))^2\; ds \)^2\right]=0.
\eeq
\ep

\bpr 
 Due to the fact that
\beq
{\cal Z}^{(n,k-1,0)}_{i,j}(\Phi,\eta-\delta_i) = {\cal Z}^{(n,k-1,0)}_{i,j}(\Phi,\eta) 
\eeq
if we  isolate the term $s=1$  in \eqref{sum} we obtain
\beq
\nabla^{i,j}\caX^{(n,k)} (\Phi,\eta)=-\cs n\(\phi(\tfrac jn) -\phi(\tfrac in)\) \cdot {\cal Z}^{(n,k-1,0)}_{i,j}(\Phi,\eta) + \widehat { G}_{i,j}^{(n,k)} (\Phi,\eta)
\eeq
then the statement follows from Proposition \ref{Prop1} and Proposition  \ref{Prop2}.
\epr

\subsubsection{Conclusion}
From \eqref{Ga} and \eqref{error} we have
\beq\label{Ga1}
&&n^2 \Gamma \caX^{(n,k)}(\Phi,\eta) = \frac{1}{n^d} \sum_{x \in \Zd } \sum_{r \in \caR} p(r) \eta_x ( \alpha + \sigma \eta_{x+r}) \( \cs\, \langle r,\nabla\phi(\tfrac{x}n)\rangle \cdot \caX^{(n,k-1)}(\Phi,\eta) - G_{x,x+r}^{(n,k)}(\phi,\eta)\)^2 \nn\\
&& = \frac{ \cs^2}{n^d}  \(\caX^{(n,k-1)}(\Phi,\eta)\)^2\cdot \sum_{x \in \Zd } \sum_{r \in \caR} |\langle r, \nabla \phi(\tfrac{x}n)\rangle|^2\;  p(r) \eta_x ( \alpha + \sigma \eta_{x+r})+ \caG_1^{(n,k)}(\Phi,\eta)\nn
\eeq
with
\beq
&& \caG_1^{(n,k)}(\Phi,\eta):=\frac{1}{n^d} \sum_{x \in \Zd } \sum_{r \in \caR} p(r) \eta_x ( \alpha + \sigma \eta_{x+r})  \cdot   G_{x,x+r}^{(n,k)}(\phi,\eta)\nn\\
 &&\hskip4cm \cdot \(G_{x,x+r}^{(n,k)}(\phi,\eta)-2 \cs \langle r,\nabla\phi(\tfrac{x}n)\rangle \caX^{(n,k-1)}(\Phi,\eta)\) \nn
\eeq
then we can write
\beq\label{Ga2}
&&n^2 \Gamma \caX^{(n,k)}(\Phi,\eta) =\\\nn
&&  = \rho(\alpha+\sigma \rho)\,\frac{ \cs^2}{n^d}  \(\caX^{(n,k-1)}(\Phi,\eta)\)^2\cdot \sum_{x \in \Zd } \sum_{r \in \caR} |\langle r, \nabla \phi(\tfrac{x}n)\rangle|^2\;  p(r)+ \caG_1^{(n,k)}(\Phi,\eta)+ \caG_2^{(n,k)}(\phi,\eta) \nn
%&&\hskip-2cm= k^2\cs^2\(\caX^{(n,k-1)}(\Phi,\eta)\)^2 \left\{  \frac { \alpha {\chi} }{n}\sum_{i \in \Z }  \eta_i   \(\phi\myprime(\tfrac{i}n)\)^2 +   \frac \sigma{n}\sum_{i \in \Z } \sum_{r \in \caR} p(r) {r^2 }\eta_i \eta_{i+r}  \(\nabla_n \phi(\tfrac{i}n)\)^2\right\}\nn
%&+&\col{\text{check-errors: There was an $r$ missing}}
\eeq
with 
\beq\label{err}
&& \caG_2^{(n,k)}(\phi,\eta) := \\
&&\frac{ \cs^2}{n^d}  \(\caX^{(n,k-1)}(\Phi,\eta)\)^2\cdot \sum_{x \in \Zd } \sum_{r \in \caR} |\langle r, \nabla \phi(\tfrac{x}n)\rangle|^2\;  p(r)\left\{\alpha(\eta_x-\rho) + \sigma (\eta_x\eta_{x+r}-\rho^2)\right\}.\nn
\eeq
We first estimate the term due to the error $ \caG_1^{(n,k)}(\Phi,\eta)$.

\bp\label{uno}
For every  $t>0$ and every test function $\phi\in S(\R)$  there exists $C>0$ such that, for all $n\in \N$,
\be\label{errbound}
 \lim_{ n \to \infty} \E_{n} \left[ \( \int_0^t  \caG_1^{(n,k)}(\Phi,\eta(n^2 s))   ds  \)^2 \right] = 0
\ee
\ep
\bpr
It follows from Proposition \ref{Prop2} and the convergence, by inductive hypothesis, of $\caX^{(n,k-1)}(\phi,\eta)$.
\epr

\noindent
The two following propositions allow us to estimate the error $ \caG_2^{(n,k)}$ and then to perform the replacement in  \eqref{Ga2}.
\bl\label{BGCarre2k1}
For every  $t>0$ and every test function $\phi\in S(\R)$  there exists $C>0$ such that, for all $n\in \N$,
\be\label{eqBGCarre2k1}
  \lim_{ n \to \infty} \E_{n} \left[ \( \int_0^t \frac{1}{n^d}\sum_{x \in \Zd } \sum_{r \in \caR} |\langle r, \nabla \phi(\tfrac{x}n)\rangle|^2\;  p(r)  (\eta_x(n^2s) -\rho)   ds  \)^2 \right] = 0.
\ee
\el
\bpr
From \eqref{one} we can write the integrand in \eqref{eqBGCarre2k1} as
\beq\label{ProofBGCarre2k1}
\frac{1}{n^{d/2}} \caX_s^{(n,1)}(\Psi), \qquad \text{with} \qquad
\Psi(\xi):=\prod_{x\in \Zd} \psi(x)^{\xi_x} \qquad \psi(x):= \sum_{r \in \caR} |\langle r, \nabla \phi(x)\rangle|^2\;  p(r) 
\eeq
then the statement follows  {from the convergence of $\caX_s^{(n,1)}(\Psi)$} and the extra factor $ \frac{1}{n^{d/2}}$.
\epr
\noindent
%From Proposition \ref{BGCarre2k1} we have
%\beq\label{exactgamma2k4}
%n^2\Gamma \caX^{(n,2)}(\phi^{(2)},\eta) 
%&=&  4\cs^2 \alpha {\chi} \rho\cdot  \(\caX^{(n,1)}(\phi,\eta)\)^2 \cdot \frac 1 n\sum_{i \in \Z }  \(\nabla_n\phi\(\tfrac i n\)\)^2 \nn \\
%&+&  {4\cs^2 \sigma }\cdot  \(\caX^{(n,1)}(\phi,\eta)\)^2 \cdot \frac 1{n}\sum_{i \in \Z } \sum_{r \in \caR} p(r) {r^2} \eta_i \eta_{i+r}  \(\nabla_n\phi\(\tfrac i n\)\)^2 \nn\\
%&+&\col{check-errors}
%\eeq
\vskip.2cm 
\noindent
Similarly, another replacement is necessary  on the second term of the RHS of \eqref{err}.
\bl\label{BGCarre2k2}
For every  $t>0$ and every test function $\phi$  we have
\be\label{EQBGCarre2k2}
\lim_{n \to \infty} \E_{n} \left[ \( \int_0^t \frac{1}{n^d}\sum_{x \in \Zd } \sum_{r\in \caR} |\langle r, \nabla \phi(\tfrac{x}n)\rangle|^2\;  p(r)\( \eta_x(n^2s) -\rho \) \( \eta_{x+r}(n^2s) -\rho\)    ds  \)^2 \right] = 0.
\ee
\el
\bpr
The proof of this lemma is done in the same spirit than Proposition \ref{metarep}.
\epr
\bp\label{due}
For every  $t>0$ and every test function $\phi\in S(\R)$  there exists $C>0$ such that, for all $n\in \N$,
\be\label{errbound}
 \lim_{ n \to \infty} \E_{n} \left[ \( \int_0^t  \caG_2^{(n,k)}(\Phi,\eta(n^2 s))   ds  \)^2 \right] = 0.
 \ee
\ep
\bpr
It follows from Lemma \ref{BGCarre2k1}, Lemma \ref{BGCarre2k2} and the convergence, by inductive hypothesis, of $\caX^{(n,k-1)}(\phi,\eta)$.
\epr
\vskip.3cm
\noindent
From Propositions  \ref{uno}  and \ref{due} we can write
\beq\label{Ga2}
&&n^2 \Gamma \caX^{(n,k)}(\Phi,\eta) =\\\nn
&&  = \rho(\alpha+\sigma \rho)\,\frac{\cs^2}{n^d}  \(\caX^{(n,k-1)}(\Phi,\eta)\)^2\cdot \sum_{x \in \Zd } \sum_{r \in \caR} |\langle r, \nabla \phi(\tfrac{x}n)\rangle|^2\;  p(r)+ \caG^{(n,k)}(\Phi,\eta)\nn
%&&\hskip-2cm= k^2\cs^2\(\caX^{(n,k-1)}(\Phi,\eta)\)^2 \left\{  \frac { \alpha {\chi} }{n}\sum_{i \in \Z }  \eta_i   \(\phi\myprime(\tfrac{i}n)\)^2 +   \frac \sigma{n}\sum_{i \in \Z } \sum_{r \in \caR} p(r) {r^2 }\eta_i \eta_{i+r}  \(\nabla_n \phi(\tfrac{i}n)\)^2\right\}\nn
%&+&\col{\text{check-errors: There was an $r$ missing}}
\eeq
where the term  $ \caG^{(n,k)}(\Phi,\eta)$ is a vanishing error:
\be
\lim_{n \to \infty} \E_{n} \left[ \( \int_0^t  \caG^{(n,k)}(\Phi,\eta(n^2s)) ds  \)^2\right]=0.\nn
\ee
Therefore we conclude that the proposed (remember that at this point we do not know if the limiting object is indeed a martingale) predictable quadratic variation of our limiting martingale is given by
\be\label{QuadvarmartingXkInduction}
 \cs^2\chi \rho (\alpha + \sigma \rho) t  \( \caY_s^{(k-1)}( \phi^{(k-1)}) \)^2 \int_{\R^d} \norm{\nabla \phi(x)}^2  dx.
\ee
\noindent
Arrived at this point we can conclude that if $\{ M_t^{(n,k)}(\Phi): t \in [0,T] \}$ has a limit as $n \to \infty$, and if the limit is a square-integrable martingale then its quadratic variation is given by \eqref{QuadvarmartingXkInduction}. In what follows we will show tightness and uniform integrability, i.e. we will prove  that $\{ M_t^{(n,k)}(\Phi): t \in [0,T] \}$ converges to $\{ M_t^{(k)}(\Phi): t \in [0,T] \}$ and that $\{ M_t^{(k)}(\Phi): t \in [0,T] \}$ is indeed a martingale.

\subsection{Tightness}\label{T}

In this section we prove tightness for the family of laws $\{ Q_n^{(k)} \}_{ n \in \N}$, induced by  $\{ \caY^{(n,k)}(\cdot,t) \}_{t \geq 0}$ on $D([0,\infty), S\myprime(\R^k))$. From  the Dynkin formula we know that
\be\label{Mtrewrit}
M_n\myprime (t,\phi^{(k)}) = \caY^{(n,k)}_t(\phi^{(k)})  - n^2 \int_0^t \loc \caY_s^{(n,k)}(\phi^{(k)}) ds
\ee
and
\be\label{Ntrew}
N_n\myprime (t,\phi^{(k)}) = M_n\myprime (t,\phi^{(k)})^2 - n^2 \int_0^t \Gamma \caY_s^{(n,k)}(\phi^{(k)}) ds
\ee
are martingales.
Theorem 2.3 in  \cite{ferrari1988non}, which we include in Appendix \ref{tightFerr}, allows us to reduce the proof the tightness of $\{ Q_n^{(k)} \}_{ n \in \N}$ to the verification  of conditions \eqref{CI}-\eqref{CII}. We verify  these conditions in Proposition \ref{Propgamma1control}, Proposition \ref{Propgamma2control} and Proposition \ref{PropModuluscontinuity2} below.

\subsubsection{The $\gamma_1$ term}
The following Proposition shows that conditions \eqref{CI} and \eqref{CIprime} hold true.
\bp\label{Propgamma1control}
For any $\phi^{(k)} \in S(\R^{kd})$ and $t_0 \geq 0$ we have:
\be\label{gamma1controlzero}
\sup_{ n \in \N} \, \sup_{ 0 \leq t \leq t_0} \E_{n} \left[ \( \caY_t^{(n,k)}(\phi^{(k)}) \)^2 \right] < \infty
\ee
and
\be\label{gamma1control}
\sup_{ n \in \N} \, \sup_{ 0 \leq t \leq t_0} \E_{n} \left[ \(n^2  \loc \caY_t^{(n,k)}(\phi^{(k)}) \)^2 \right] < \infty.
\ee
\ep
\bpr
We start with the proof  \eqref{gamma1control} which is more involved. Thanks to stationarity, the expectation does not depend on time, and then  we can ignore the supremum over time in \eqref{gamma1control}.
From   \eqref{closeddrifttwo} we already have an expression for the integrand of \eqref{gamma1control}:
\be\label{aftermathgamma1}
n^2 \loc \caY^{(n,k)}(\phi^{(k)},\eta) = \alpha k \cdot\tfrac {\chi}2 \cdot \caY^{(n,k)}(\phi^{(k-1)}\otimes \Delta\phi,\eta) + O(n^{-1})
\ee
recall that here again we are using the fact that the field $\caY^{(n,k)}$ can be also thought as acting on general ( not necessarily symmetric ) test functions.
Because of stationarity it is enough to estimate
\be
\E_{\nu_\rho} \left[ \( \caY^{(n,k)}(\phi^{(k-1)}\otimes \Delta\phi,\eta) \)^2 \right].
\ee
then the  desired bound is obtained by applying Proposition \ref{mombndXk}.
In the same spirit we can use Proposition \ref{mombndXk}   to bound \eqref{gamma1controlzero}. 
\epr

\subsubsection{The $\gamma_2$ term}
 Similarly to the previous section, here we prove the  following proposition in order to verify the condition \eqref{CIprime} for $\gamma_2$.
\bp\label{Propgamma2control}
For any $\phi^{(k)} \in S(\R^{kd})$ and $t_0 \geq 0$ we have:
\be\label{gamma2control}
\sup_{ n \in \N} \, \sup_{ 0 \leq t \leq t_0} \E_{n} \left[ \(n^2 \Gamma \caY_t^{(n,k)}(\phi^{(k)}) \)^2 \right] < \infty.
\ee
\ep
\bpr
Thanks to stationarity we can neglect the supremum over time. Recall that in \eqref{Ga2} we have obtained  an  expression for the integrand on \eqref{gamma2control}
\beq\label{gamma2control2}
n^2\Gamma \caY^{(n,k)}(\phi^{(k)},\eta)&=& \nn \\
&&\rho(\alpha+\sigma \rho)\,\frac{k^2 \cs^2}{n^d}  \(\caX^{(n,k-1)}(\Phi,\eta)\)^2\cdot \sum_{x \in \Zd } \sum_{r \in \caR} |\langle r, \nabla \phi(\tfrac{x}n)\rangle|^2\;  p(r) + O(n^{-1}) \nn 
\eeq
taking the square of which we obtain
\beq\label{Gamma2control3}
&& \E_{\nu_\rho} \left[ (n^2\Gamma \caY^{(n,k)}(\phi^{(k)},\eta))^2 \right] \nn \\
&=& \rho^2(\alpha+\sigma \rho)^2\,\frac{k^4 \cs^4}{n^{2d}}  \cdot \sum_{x, y \in \Zd } \sum_{r_1, r_2 \in \caR} |\langle r_1, \nabla \phi(\tfrac{x}n)\rangle|^2 \cdot |\langle r_2, \nabla \phi(\tfrac{y}n)\rangle|^2\;  p(r_1) \cdot p(r_2) \nn \\
&\times& \E_{\nu_\rho} \left[ \(\caX^{(n,k-1)}(\Phi,\eta)\)^4 \right].
\eeq
Notice that the first factor on the RHS of \eqref{Gamma2control3} can be controlled by using the compact support of $\phi$ and the factor $\frac{1}{n^{2d}}$. It is then sufficient to estimate 
\be
\sup_{ n \in \N} \, \E_{\nu_\rho} \left[ \(\caX^{(n,k-1)}(\Phi,\eta)\)^4 \right]
\ee
then Proposition \ref{mombndXk} finishes the proof.
\epr

\subsubsection{Modulus of continuity}
In this section we show that condition (2.5) of Theorem 2.3 in  \cite{ferrari1988non} is satisfied.
\bp\label{PropModuluscontinuity2}
For every $\phi^{(k)} \in S(\R^{kd})$ there exists a sequence $\delta(t,\phi,n)$ converging to zero as $n \to 0$ such that:
\be
\lim_{ n \to \infty} \Pr_{n} ( \sup_{ 0 \leq t \leq T} \lvert \caY_t^{(n,k)}(\phi^{(k)},\eta) -  \caY_{t-}^{(n,k)}(\phi^{(k)},\eta)  \rvert \geq \delta(t,\phi,n)) = 0.
\ee 
\ep
\bpr
Let us choose a sequence $\{ \delta_n \}_{n \geq 1}$, independent of $\phi$ and $t$, and converging to zero as $n \to \infty$.   For the moment we do not specify the rate of convergence of the sequence.\\
\noindent
By the Dynkin formula we have
\be
\sup_{ 0 \leq t \leq T} \lvert \caY_t^{(n,k)}(\phi^{(k)}) -  \caY_{t-}^{(n,k)}(\phi^{(k)})  \rvert = \sup_{ 0 \leq t \leq T} \lvert M_t^{(n,k)}(\phi^{(k)}) -M_{t-}^{(n,k)}(\phi^{(k)}) \rvert 
\ee
which, together with Doobs's inequality for sub-martingales, implies
\beq
&& \Pr_{n} \left( \sup_{ 0 \leq t \leq T} \lvert \caY_t^{(n,k)}(\phi^{(k)}) -  \caY_{t-}^{(n,k)}(\phi^{(k)})  \rvert \geq \delta_n \right) \nn \\
&\leq& \frac{1}{\delta_n^2} \, \E_{n} \(  \( M_T^{(n,k)}(\phi^{(k)}) -M_{T-}^{(n,k)}(\phi^{(k)})   \)^2 \) \nn \\
&\leq& \frac{1}{\delta_n^2} \, \E_{n} \(   \int_{T-}^{T} n^2 \Gamma \caY_s^{(n,k)}(\phi^{(k)}) ds \) \nn \\
&\leq& \frac{1}{\delta_n^2} \, \E_{n} \(   \int_{T-\epsilon_n}^{T} n^2 \Gamma \caY_s^{(n,k)}(\phi^{(k)}) ds \) \nn \\
&\leq& \frac{\epsilon_n}{\delta_n^2} C(\phi,\rho,\alpha)
\eeq
where in the third inequality we used the positivity of the carr\'e-du-champ, and in the last inequality we used the stationarity of the measure, together with an  estimate similar to the one used to show Proposition \ref{Propgamma2control}.
Thus, to conclude the proof, it is sufficient to choose $\epsilon_n$ such that
\be\label{epsilonnchoice}
\lim_{ n \to \infty} \frac{\epsilon_n}{\delta_n^2} = 0.
\ee
\epr	
\br
Notice that  Proposition \ref{PropModuluscontinuity2} implies,  in particular, that the law induced by $\caX_t^{(n,k)}$ is concentrated in continuous paths.
\er

\subsection{Characterization of limit points}\label{Char}

At this point we can only say that the sequence $\{ M_t^{(n,k)}(\cdot): t \in [0,T] \}$ converges weakly to the process $\{ M_t^{(k)}(\cdot): t \in [0,T] \}$ satisfying expressions \eqref{martingXk} and \eqref{QuadvarmartingXk}. Nevertheless, we would like to support the claim, given in Theorem \ref{mainkthorder}, that the limiting process $\{ M_t^{(k)}(\cdot): t \in [0,T] \}$ is indeed  a martingale with the proposed predictable quadratic variation given by \eqref{QuadvarmartingXkInduction}. At this aim  we prove the following result.

\bp\label{UIofMn}
The sequence $\{ M_t^{(n,k)}(\cdot): t \in [0,T] \}$ is uniformly integrable.
\ep

\bpr
By standard arguments it is enough to provide a uniform $L^p(\P_n)$ bound for $p >1$. Notice that, thanks to the martingale decomposition \eqref{Mtdef},  and the same type of arguments used in the proofs of Propositions \ref{Propgamma1control} and \ref{Propgamma2control}, we can indeed find the desired bounds for $p=2$.
\epr

%\col{Here, when saying standard arguments we mean the following:
%\ben
%\item We already have convergence of distributions
%\item If necesseary we move to subsequences to make things match
%\item By Skohorods we have almost sure convergence of copies of the $M_n$ in some abstract probability space
%\item Vitalis theorem (UI) gives convergence in mean of the copies of $M_n$ to $M$
%\item From convergence in mean we obtain convergence of the means.
%\een}

\noindent
The same type of reasoning used in Proposition \ref{UIofMn} gives us the following result.

\bp\label{UIofNn}
The sequence $\{ N_t^{(n,k)}(\cdot): t \in [0,T] \}$ is uniformly integrable.
\ep

\noindent
Combining Propositions \ref{UIofMn} and \ref{UIofNn} we show that any limit point of the sequence $\{ M_t^{(n,k)}(\cdot): t \in [0,T] \}$ satisfies the recursive martingale problem \eqref{martingXk}-\eqref{QuadvarmartingXk}. 

\noindent
\subsection{Uniqueness}\label{Uni}
It remains to show uniqueness of the solution of the martingale problem \eqref{martingXk}-\eqref{QuadvarmartingXk}. First notice that by Duhamel formula, from \eqref{martingXk}, we can deduce
\be\label{Duhamels}
\caY_t^{(k)}(\phi^{(k)}) = \caY_0^{(k)}(S_t \phi^{(k)}) + \int_0^t d M_s^{(k)}(S_{t-s} \phi^{(k)})
\ee
where $S_t$ is the semigroup associated to the $kd$ dimensional Laplacian. Given the distribution of $\caY_0^{(k)}$ and the well-definiteness of $M_t^{(k)}$, the RHS of \eqref{Duhamels} uniquely determines  the  finite dimensional distributions of $\caY_t^{(k)}$.  Then, by the continuity of $\caY_t^{(k)}$, we conclude the uniqueness of limiting point.

\section{Appendix}

\subsection{Carr\'e-du-champ}

\bp\label{squareduchamp}
Consider an interacting  particles system with generator 
\be
L f (\eta) = \sum_{\eta\myprime} c(\eta, \eta\myprime) \( f(\eta\myprime) - f(\eta) \)
\ee
the following is an alternative formulation for its carr\'e-du-champ
\be\label{Eqsquareduchamp}
\Gamma (f)(\eta) = \sum_{\eta\myprime} c(\eta, \eta\myprime) \( f(\eta\myprime) - f(\eta) \)^2.
\ee
\ep

\bpr
By definition we have 
\beq
\Gamma (f)(\eta) &=& \sum_{\eta\myprime} c(\eta, \eta\myprime) \( f(\eta\myprime)^2 - f(\eta)^2 \) \nn \\
&-& 2 f(\eta) \sum_{\eta\myprime} c(\eta, \eta\myprime) \( f(\eta\myprime) - f(\eta) \) \nn \\
&=& \sum_{\eta\myprime} c(\eta, \eta\myprime) \( f(\eta\myprime)^2 - f(\eta)^2 \) \nn \\
&-& \sum_{\eta\myprime} c(\eta, \eta\myprime) \( 2 f(\eta)f(\eta\myprime) - 2f(\eta)^2 \) \nn \\
&=& \sum_{\eta\myprime} c(\eta, \eta\myprime) \( f(\eta\myprime)^2 - f(\eta)^2 -2 f(\eta)f(\eta\myprime) + 2f(\eta)^2\) \nn \\
&=& \sum_{\eta\myprime} c(\eta, \eta\myprime) \( f(\eta\myprime)^2  -2 f(\eta)f(\eta\myprime) + f(\eta)^2\) \nn \\
&=& \sum_{\eta\myprime} c(\eta, \eta\myprime) \( f(\eta\myprime)  - f(\eta)\)^2 \nn
\eeq
that concludes the proof.
\epr

\subsection{Tightness criterium}\label{tightFerr}

In this section we state a well known criterium for tightness extracted from \cite{ferrari1988non}:

\bt
Let $(\Omega, \caF)$ be a measurable space with right-continuous filtrations $\{ \caF_t^{n} \}_{t \geq 0}$ and probability measures $\Pr_{n}(\cdot)$, $ n \in \N$. Let $\{ \caX_t^{n} \}_{t \geq 0}$ be an $\caF_t^{n}$-adapted process with paths in $D([0,\infty), S\myprime(\R^k))$ and let us also suppose that there exists, 
for each $\phi \in S(\R^k)$, $\caF_t^{n}$-predictable processes $\gamma_1^{n}(\cdot,\phi), \gamma_2^{n}(\cdot,\phi)$ such that:
\be\label{Mttightness}
M_t^{n} (\phi) := \caX_t^{n}(\phi) -  \int_0^t \gamma_1^{n}(s,\phi) ds
\ee
and
\be\label{Nttightness}
M_t^{n} (\phi)^2 -  \int_0^t \gamma_2^{n}(s,\phi) ds
\ee
are martingales. Assume further that  it holds:
\begin{description}
\item[CI:] for $t_0 \geq 0$ and $\phi \in S(\R^k)$:
\be\label{CI}
\sup_{ n \in \N} \, \sup_{ 0 \leq t \leq t_0} \E_{n} ( \caX_t^{n}(\phi)^2) < \infty
\ee
and for $ i \in \{1,2\}$:
\be\label{CIprime}
\sup_{ n \in \N} \, \sup_{ 0 \leq t \leq t_0} \E_{n} ( \gamma_i^{n}(t,\phi)^2) < \infty;
\ee
\item[CII:] for every $\phi \in S(\R^k)$ there exists a sequence $\delta(t,\phi,n)$ converging to zero as $n \to 0$ such that:
\be\label{CII}
\lim_{ n \to \infty} \Pr_{n} ( \sup_{ 0 \leq s \leq t} \lvert \caX_s^{n}(\phi) - \caX_{s-}^{n}(\phi) \rvert \geq \delta(t,\phi,n)) = 0
\ee 
\end{description}
then the family of laws $\{ Q^{n} \}_{ n \in \N}$, induced by  $\{ \caX_t^{n} \}_{t \geq 0}$ on $D([0,\infty), S\myprime(\R^k))$ under $\P_n$, is a tight family and any weak limit point is supported by $C([0,\infty), S\myprime(\R^k))$.
\et

\section*{Acknowledgements}
The authors would like to thank Federico Sau for helpful discussions; M. Ayala acknowledges financial support from the Mexican Council on Science and Technology (CONACYT) via the scholarship 457347.

%\nocite{*}
\bibliographystyle{abbrv}
\bibliography{../Sticky/Bibliography/Biblio}

\end{document}